\documentclass[10pt]{article}
\usepackage{amsmath, amssymb, amscd, amsthm}

\newcommand{\xra}[1]{\ensuremath{\xrightarrow{#1}}}
\newcommand{\B}[1]{\ensuremath{\mathbb{#1}}}
\newcommand{\C}[1]{\ensuremath{\mathcal{#1}}}
\newcommand{\G}[1]{\ensuremath{\mathfrak{#1}}}
\newcommand{\eotimes}[1]{\ensuremath{\underset{#1}{\otimes}}}
\newcommand{\fotimes}[1]{\ensuremath{\underset{#1}{\Box}}}

\newtheorem{thm}{Theorem}[section]
\newtheorem{cor}[thm]{Corollary}
\newtheorem{lem}[thm]{Lemma}

\theoremstyle{definition}
\newtheorem{defn}[thm]{Definition}
\newtheorem{rem}[thm]{Remark}

\numberwithin{equation}{section}

\title{Bialgebra Cyclic Homology with Coefficients\\ Part I}
\author{Atabey Kaygun}
\date{}
\begin{document}
\maketitle

\section{Introduction}

Cyclic cohomology of Hopf algebras admitting a modular pair was first
defined in \cite{ConnesMoscovici:HopfCyclicCohomology} and further
developed in \cite{ConnesMoscovici:HopfCyclicCohomologyI} and
\cite{ConnesMoscovici:CyclicCohomologyOfHopfAlgebras} in the context of
transverse geometry.  Their results are followed by several papers
computing Hopf cyclic (co)homology of certain Hopf algebras such as
\cite{Taillefer:CyclicHomologyOfHopfAlgebras},
\cite{Crainic:CyclicCohomologyOfHopfAlgebras} and
\cite{KhalkhaliRangipour:CyclicCohomologyOfExtendedHopf}.  In a series
of papers \cite{KhalkhaliAkbarpour:CyclicHomologyHopfModuleAlgebras},
\cite{Khalkhali:InvariantCyclicHomology},
\cite{Khalkhali:DualCyclicHomology}, and
\cite{Khalkhali:HopfCyclicHomology} the authors developed a theory of
cyclic (co)homology which works with Hopf module coalgebras or Hopf
comodule algebras with coefficients in stable Hopf module/comodules
satisfying anti-Yetter-Drinfeld condition (SaYD.) In this paper, We show
that one can extend the Hopf cyclic homology non-trivially by using just
stable module/comodules, dropping the aYD condition.  This also allows
us to extend the definition of the cyclic homology to bialgebras.

The version of the Connes--Moscovici cyclic homology presented in
Section~\ref{OldHomology} and the extension we provided in
Section~\ref{NewHomology} use only module coalgebras and their module
coinvariants over a bialgebra $B$.  One can use comodule algebras and
comodule invariants over a bialgebra to construct a dual theory as in
\cite{Khalkhali:DualCyclicHomology}.  In the case of Hopf algebras and
SaYD coefficient modules, these theories are dual to each other in the
sense of \cite{Khalkhali:CyclicDuality}.  We consider this dual theory
in a separate paper \cite{Kaygun:BialgebraCyclicII}.

Here is a plan of this paper: In Section~\ref{convention}, we set up the
notation and list overall assumptions we make.  In
Section~\ref{OldHomology}, we give a self-contained account of the
Connes--Moscovici cyclic homology, as it is developed by Khalkhali
et. al.  The presentation basically follows the papers mentioned above
in results, but differs in proofs and other details.  However, some of
the results in Section~\ref{OldHomology} are new, such as
Theorem~\ref{Coadjoint}.  In Section~\ref{NewHomology}, we define the
bialgebra cyclic homology with stable coefficients.
Theorem~\ref{main_isomorphism}, and discussion following it, allows us
to define the new bialgebra cyclic homology.
Theorem~\ref{generalization} shows the relationship between the new
bialgebra cyclic homology and the classical Hopf cyclic homology.
Section~\ref{TechnicalResults} contains the necessary technical results
needed for computations.  Section~\ref{Computations} contains
computations of the bialgebra cyclic homology of two Hopf algebras with
stable but non-aYD coefficients.  In
Section~\ref{TransverseCodimensionN} there are two computations
involving Hopf algebra of transverse geometry of foliations in
codimension $n$, namely $\C{H}(n)$, with two different coefficient
modules.  Section~\ref{quantum} contains a computation of the bialgebra
cyclic homology of $U_q(\G{g})$, quantum deformation of an arbitrary
semi-simple Lie algebra $\G{g}$, with coefficient in an arbitrary copy
of $U_qsl(2)$ embedded in $U_q(\G{g})$.

\section{Notation and conventions}\label{convention}

\noindent We assume $k$ is a field of arbitrary characteristic and $H$
is a Hopf algebra over $k$.  In the case of a regular commutative
algebra $k$, most of the results of this paper can be obtained if $H$ is
an injective and a flat $k$--module.  We also assume $H$ has a bijective
antipode.

\noindent Whenever we refer an object ``simplicial'' or
``cosimplicial,'' the reader should read as ``pre-simplicial'' and
``pre-cosimplicial'' meaning that we did not consider (co)degeneracy
morphisms as a part of the (co)simplicial data.

\noindent A simplicial $X_*$ module is called para-(co)cyclic iff it is
almost a (co)cyclic module, in that it satisfies all axioms of a
(co)cyclic module except that the action of $\tau_n$ on each $X_n$ need
not to be of order $n+1$, for any $n\geq 0$.

\noindent A (para-)(co)cyclic module $\C{Z}_*$ is called a
(para-)(co)cyclic $H$--module iff all structure morphisms are
$H$--module morphisms.

\noindent The tensor product over an algebra $A$ is denoted by
$\eotimes{A}$, and a cotensor product over a coalgebra $C$ is denoted
by $\fotimes{C}$.  Recall that if $X\xra{\rho_X}X\otimes C$ and
$Y\xra{\rho_Y}C\otimes Y$ are two $C$--comodules (right and left
respectively), then $X\fotimes{C}Y$ is defined as $ker((\rho_X\otimes
id_Y)-(id_X\otimes \rho_Y))$.  The derived functors of the cotensor
product are denoted by $\text{Cotor}_C^*(\cdot,\cdot)$.

\noindent Given a counital bialgebra $(B,\cdot,\Delta,\epsilon)$ and a
left $B$--module $M$, the module of $B$--coinvariants of $M$ which is
\[ M/\left\{\epsilon(b)m-b\cdot m|\ b\in B,\ m\in M\right\} \] 
is denoted by ${}_BM$.

\noindent For a coalgebra $(C,\Delta)$, we use Sweedler's notation and
denote $\Delta(c)$ by $c_{(1)}\otimes c_{(2)}$.  Similarly, for a (left)
$C$--comodule $Y\xra{\rho_Y}C\otimes Y$, we use $\rho_Y(y) =
y_{(-1)}\otimes y_{(0)}$ for the coaction morphism.

\section{The Connes--Moscovici cyclic homology}\label{OldHomology}

\begin{defn}
Let $H$ be a Hopf algebra and $Y$ be a $H$--module/comodule.  Consider
the graded $k$--module $\left\{H^{\otimes n+1}\otimes Y\right\}_{n\geq
0}$ with the following $k$--module homomorphisms
\begin{align*}
\partial_0(h^0\otimes\cdots\otimes h^n\otimes y)
 = & (h^0_{(1)}\otimes h^0_{(2)}\otimes\cdots\otimes y)
\end{align*}
Moreover, we define
\begin{align*}
\tau_n(h^0\otimes\cdots\otimes h^n\otimes y)
 = \left(S^{-1}(y_{(-1)})h^n\otimes h^0\otimes\cdots\otimes 
         h^{n-1}\otimes y_{(0)}\right)
\end{align*}
It is easy to check that
\begin{align*}
\tau_n^{-1}(h^0\otimes\cdots\otimes h^n\otimes y)
 = (h^1\otimes\cdots\otimes h^n\otimes y_{(-1)}h^0\otimes y_{(0)})
\end{align*}
is the inverse of $\tau_n$ for any $n\geq 1$.  If we define
\begin{align}\label{cosimplicial_structure}
\partial_j
 = & \tau_{n+1}^j\partial_0\tau_n^{-j}
\end{align}
for $0\leq j\leq n+1$, one can easily check that
\begin{align*}
\partial_i\partial_j =& \partial_{j+1}\partial_i 
	\hspace{1cm}\text{ if } i\leq j
\end{align*}
which means we have a cosimplicial object in the category of
$k$--modules.  On the other hand, for $0\leq j\leq n$
\begin{align*}
\tau_{n+1}\partial_j & = \partial_{j+1}\tau_n
\end{align*}
and for $j=n+1$,
\begin{align*}
\tau_{n+1}\partial_{n+1} 
& = \tau_{n+1}^{n+2}\partial_0\tau_n^{-n-1}=\partial_0
\end{align*}
since 
\begin{align*}
\partial_0\tau_n^{-n-1}
  (h^0\otimes\cdots\otimes h^n\otimes y) 
 = & \partial_0\left(y_{(-n-1)}h^0\otimes\cdots\otimes
         y_{(-1)}h^n\otimes y_{(0)}\right)\\
 = & \left(y_{(-n-1)(1)}h^0_{(1)}\otimes 
         y_{(-n-1)(2)}h^0_{(2)}\otimes\cdots\otimes y_{(0)}\right)\\
 = & \tau_{n+1}^{-n-2}(h^0_{(1)}\otimes h^0_{(2)}\otimes\cdots\otimes y)\\
 = & \tau_{n+1}^{-n-2}\partial_0(h^0\otimes\cdots\otimes h^n\otimes y)
\end{align*}
This means, the collection $\left\{H^{\otimes n+1}\otimes
Y\right\}_{n\geq 0}$ is a para-cocyclic module since $\tau_n^{n+1}$
fails to be $id_n$.  We denote this object by $\B{T}_*(H,Y)$.
\end{defn}

\begin{defn}
Given a coalgebra $C$ and a left $C$--comodule $X$ and a right
$C$--comodule $Y$, one can form the cobar cosimplicial complex
$B_*(X,C,Y)=\{X\otimes C^{\otimes n}\otimes Y\}_{n\geq 0}$ where the
cosimplicial structure morphisms are
\begin{align*}
d_j(x\otimes c^1\otimes\cdots\otimes c^n\otimes y)
 = \begin{cases}
   \left(x_{(0)}\otimes x_{(1)}\otimes c^1\otimes\cdots\otimes c_n
         \otimes y\right)        & \text{ if } j=0\\
   \left(x\otimes\cdots\otimes c^j_{(1)}\otimes c^j_{(2)}\otimes\cdots
         \otimes y\right)        & \text{ if } 0<j\leq n\\
   \left(x\otimes c^1\otimes\cdots\otimes c^n\otimes y_{(-1)}
         \otimes y_{(0)}\right)  & \text{ if } j=n+1
   \end{cases}
\end{align*}
\end{defn}

\begin{defn}
For any $n\geq 1$, define a right action of $C$ on $C^{\otimes n}$ by
\begin{align*}
(c^1\otimes\cdots\otimes c^n)\cdot  c 
 = (c^1c_{(1)}\otimes\cdots\otimes c^n c_{(n)})
\end{align*}
for all $c\in C$ and $(c^1\otimes\cdots\otimes c^n)\in C^{\otimes n}$.
One can similarly define a left action by 
\begin{align*}
c\cdot (c^1\otimes\cdots\otimes c^n)
 = (c_{(1)}c^1\otimes\cdots\otimes  c_{(n)}c^n)
\end{align*}
We also refer left (resp. right action) of an element $c\in C$ on
$C^{\otimes n}$ by $L_c$ (resp. $R_c$) for any $n\geq 1$.
\end{defn}

\begin{lem}
Let $M$ be a $H$--bicomodule where left and right comodule structures
are denoted by
\begin{align*}
\rho_L(m) = & \left(m_{(-1)}\otimes m_{(0)}\right)\\
\rho_R(m) = & \left(m_{(0)}\otimes m_{(1)}\right)
\end{align*}
Then $M$ can also be thought as a right $H$--comodule via the coadjoint
coaction defined as
\begin{align*}
\rho_{coad}(m) 
   = & \left(m_{(0)}\otimes m_{(1)}S(m_{(-1)})\right)
\end{align*}
for any $m\in M$.
\end{lem}

\begin{thm}\label{Coadjoint}
Assume $H$ is a Hopf algebra.  Then, there is an isomorphism of
cosimplicial $k$--modules $\B{T}_*(H,Y)\xra{\Phi_*}B_*(coad(H),H,Y)$
where
\begin{align*}
\Phi_n(h^0\otimes\cdots\otimes h^n\otimes y)
 = & \left(h^0_{(2)}\otimes (h^1\otimes\cdots\otimes h^n)\cdot 
            S(h^0_{(1)})\otimes y\right)
\end{align*}
for all $(h^0\otimes\cdots\otimes h^n\otimes y)$ from $\B{T}_n(H,Y)$ and
for all $n\geq 0$.
\end{thm}
\begin{proof}
First, let me show that this is an isomorphism of graded $k$--modules:
The inverse of $\Phi_*$ is given by
\begin{align*}
\Phi_n^{-1}(h^0\otimes\cdots\otimes h^n\otimes y)
 = & \left(h^0_{(2)}\otimes (h^1\otimes\cdots\otimes h^n)
           \cdot  h^0_{(1)}\otimes y\right)
\end{align*}
for all $(h^0\otimes\cdots\otimes h^n\otimes y)$ from $B_n(H,Y)$ and for
all $n\geq 0$.  One can see this by observing
\begin{align*}
\Phi_n^{-1}\Phi_n(h^0\otimes\cdots\otimes h^n\otimes y)
  = & \Phi_n\left(h^0_{(2)}\otimes (h^1\otimes\cdots\otimes h^n)
           \cdot  S(h^0_{(1)})\otimes y\right)\\
  = & \left(h^0_{(2)(2)}\otimes (h^1\otimes\cdots\otimes h^n)\cdot 
           S(h^0_{(1)})h^0_{(2)(1)}\otimes y\right)\\
  = & (h^0\otimes\cdots\otimes h^n\otimes y)
\end{align*}
The proof that $\Phi_n\Phi_n^{-1}=id_n$ is similar.

The cosimplicial structure on $B_*(coad(H),H,Y)$ is given by
\begin{align*}
d_j(h^0\otimes h^1\otimes\cdots\otimes h^n\otimes y)
 = \begin{cases}
   (h^0_{(2)}\otimes h^0_{(3)} S(h^0_{(1)})\otimes 
    h^1\otimes\cdots\otimes h^n\otimes y)
           & \text{ if } j=0\\
   (h^0\otimes\cdots\otimes h^j_{(1)}\otimes h^j_{(2)}\otimes\cdots
     \otimes y)
           & \text{ if } 1\leq j\leq n\\
   (h^0\otimes\cdots\otimes h^n\otimes y_{(-1)}\otimes y_{(0)})
           & \text{ if } j=n+1
   \end{cases}
\end{align*}
The cosimplicial maps on $\B{T}_*(M,Y)$ for $n\geq 0$ are given by
\begin{align*}
\partial_j(h^0\otimes\cdots\otimes h^n\otimes y)
= & \begin{cases}
      (\cdots\otimes h^j_{(1)}\otimes h^j_{(2)}\otimes\cdots\otimes y)
                 & \text{ if } 0\leq j\leq n\\
      \left(h^0_{(2)}\otimes h^1\otimes\cdots\otimes h^n\otimes y_{(-1)}
            h^0_{(1)}\otimes y_{(0)}\right)
                 & \text{ if } j=n+1
    \end{cases}
\end{align*}
from the definition of cosimplicial structure morphisms given
in~(\ref{cosimplicial_structure}). Then, for $1\leq j\leq n$ we have
\begin{align*}
\Phi_{n+1}\partial_j(h^0\otimes\cdots\otimes h^n\otimes y)
 = & \Phi_{n+1}\left(h^0\otimes\cdots\otimes h^j_{(1)}\otimes
               h^j_{(2)}\otimes\cdots\otimes y\right)\\
 = & \left(h^0_{(n+2)}\otimes\cdots\otimes 
           h^j_{(1)} S(h^0_{(n+1-j)})\otimes
           h^j_{(2)} S(h^0_{(n-j)})\otimes\cdots\otimes y\right)\\
 = & d_j\Phi_n(h^0\otimes\cdots\otimes h^n\otimes y)
\end{align*}
However, for $j=0$, we have
\begin{align*} 
d_0\Phi_n (h^0\otimes\cdots\otimes h^n\otimes y)
 = & d_0\left(h^0_{(2)}\otimes (h^1\otimes\cdots\otimes
            h^n)\cdot  S(h^0_{(1)})\otimes y\right)\\
 = & \left(h^0_{(3)}\otimes h^0_{(4)}S(h^0_{(2)})\otimes 
           (h^1\otimes\cdots\otimes h^n)\cdot  S(h^0_{(1)})
            \otimes y\right)\\
 = & \left(h^0_{(2)}\otimes (h^0_{(3)}\otimes h^1\otimes\cdots\otimes 
           h^n)\cdot  S(h^0_{(1)})\otimes y\right)\\
 = & \Phi_{n+1}\partial_0(h^0\otimes\cdots\otimes h^n\otimes y)
\end{align*}
and finally for $j=n+1$,
\begin{align*}
\partial_{n+1}\Phi_n^{-1}(h^0\otimes\cdots\otimes h^n\otimes y)
 = & \partial_{n+1}\left(h^0_{(2)}\otimes(h^1\otimes\cdots\otimes h^n)
         \cdot  h^0_{(1)}\otimes y\right)\\
 = & \left(h^0_{(3)}\otimes (h^1\otimes\cdots\otimes h^n)
         \cdot  h^0_{(1)} \otimes y_{(-1)} h^0_{(2)}\otimes
          y_{(0)}\right)\\
 = & \Phi_{n+1}^{-1}\left(h^0\otimes\cdots\otimes h^n\otimes y_{(-1)}
         \otimes y_{(0)}\right)\\
 = & \Phi_{n+1}^{-1}d_{n+1}(h^0\otimes\cdots\otimes h^n\otimes y)
\end{align*}
which proves that $\Phi_*$ is an isomorphism of cosimplicial modules.
\end{proof}

\begin{rem}
The cosimplicial module $\B{T}_*(H,Y)$ has also a left $H$--module
structure.  We can transport this structure to $B_*(coad(H),H,Y)$ by
using $\Phi_*$.  The left action of $H$ on $B_*(coad(H),H,Y)$ is defined
as follows:
\begin{align*}
L_h(h^0\otimes h^1 & \otimes\cdots\otimes h^n\otimes y)\\
 := & \Phi_n L_h\Phi^{-1}_n(h^0\otimes h^1\otimes\cdots\otimes
                            h^n\otimes y) \\
  = & \Phi_n\left(h_{(1)}h^0_{(2)}\otimes h_{(2)}\cdot\left(h^1 
                  \otimes\cdots\otimes h^n\right)\cdot  h^0_{(1)}
		  \otimes h_{(3)}\cdot y\right)\\
  = & \left(h_{(1)(2)}h^0_{(2)(2)}\otimes h_{(2)}\cdot\left(h^1 
                  \otimes\cdots\otimes h^n\right)\cdot 
                   h^0_{(1)}S(h^0_{(2)(1)})S(h_{(1)(1)})
		  \otimes h_{(3)}\cdot y\right)\\
  = & \left(h_{(2)}h^0\otimes h_{(3)}\cdot\left(h^1 
                  \otimes\cdots\otimes h^n\right)\cdot 
                  S(h_{(1)})\otimes h_{(4)}\cdot y\right)
\end{align*}
for any $h\in H$ and for any $(h^0\otimes\cdots\otimes h^n\otimes y)$
from $B_*(coad(H),H,Y)$
\end{rem}

\begin{defn}
The cosimplicial module $\B{T}_*(H,Y)$ has a para-cocyclic structure,
and one can transport this structure on $B_*(coad(H),H,Y)$ by using the
isomorphism $\Phi_*$.  Then, the action of $\B{Z}$ on $B_*(coad(H),H,Y)$
is defined as
\begin{align*}
t_n^{-1} (h^0\otimes\cdots\otimes h^n\otimes y)
:= & \Phi_n\tau_n^{-1}\Phi_n^{-1}(h^0\otimes\cdots\otimes h^n\otimes y)\\
 = & \Phi_n\tau_n^{-1}\left(h^0_{(2)}\otimes (h^1\otimes\cdots\otimes h^n)
        \cdot  h^0_{(1)}\otimes y\right)\\
 = & \Phi_n\left((h^1\otimes\cdots\otimes h^n)\cdot 
        h^0_{(1)}\otimes y_{(-1)}h^0_{(2)}\otimes y_{(0)}\right)\\
 = & \Phi_n\left(h^1h^0_{(1)}\otimes(\cdots\otimes h^n
        \otimes y_{(-1)})\cdot  h^0_{(2)}\otimes y_{(0)}\right)\\
 = & \left(h^1_{(2)}h^0_{(2)}\otimes(\cdots\otimes h^n\otimes y_{(-1)})
        \cdot  h^0_{(3)}S(h^0_{(1)})S(h^1_{(1)})
        \otimes y_{(0)}\right)
\end{align*}
for any element $(h^0\otimes\cdots\otimes h^n\otimes y)$ from the
cosimplicial $k$--module $B_*(coad(H),H,Y)$.  This new para-cocyclic
module is denoted by $BC_*(coad(H),H,Y)$ and is called the para-cocyclic
cobar complex of $H$ twisted by the $H$--module/comodule $Y$.
\end{defn}

\begin{rem}
We note that $B_*(coad(H),H,Y)$ is not a cocyclic, not even a
cosimplicial, $H$--module.
\end{rem}

\begin{defn}
Now, define a graded $k$--module $\B{CM}_*(H,Y)=\left\{H^{\otimes
n}\otimes Y\right\}_{n\geq 0}$ and a morphisms of graded $k$--modules of
the form $\B{T}_*(H,Y)\xra{p_*}\B{CM}_*(H,Y)$ by
\begin{align*}
p_n(h^0\otimes\cdots\otimes h^n\otimes y)
 = \begin{cases}
   S(h^0) y & \text{ if } n=0\\
   \left(S(h^n_{(n+1)})y_{(-n)}h^0\otimes\cdots\otimes
         S(h^n_{(2)})y_{(-1)}h^{n-1}\otimes S(h^n_{(1)})y_{(0)}\right)
           & \text{ if } n>0
   \end{cases}
\end{align*}
Define also $\B{CM}_*(H,Y)\xra{i_*}\B{T}_*(H,Y)$ by
\begin{align*}
i_n(h^0\otimes\cdots\otimes h^{n-1}\otimes y)
 = (h^0\otimes\cdots\otimes h^{n-1}\otimes y_{(-1)}\otimes y_{(0)})
\end{align*}
\end{defn}

\begin{defn}
Assume $H$ is a Hopf algebra.  Then a $H$--module/comodule $Y$ is called
$m$-stable if \[ S^m(y_{(-1)})y_{(0)}=y \] for all $y\in Y$.  If $Y$ is
both $0$--stable and $1$--stable, we call it stable.
\end{defn}

\begin{defn}
Assume $H$ is a Hopf algebra.  Then $H$--module/comodule is called
anti-Yetter-Drinfeld (aYD) module iff
\begin{align*}
(hx)_{(-1)}\otimes (hx)_{(2)} = h_{(1)}x_{(-1)}S^{-1}(h_{(3)}\otimes
  h_{(2)}x_{(0)} 
\end{align*}
for any $x\in X$ and $h\in H$.  
\end{defn}

\begin{lem}
Assume $X$ is an anti-Yetter-Drinfeld module.  Then $x_{(-1)}x_{(0)}=x$
for any $x\in X$ iff $S(x_{(-1)})x_{(0)}=x$ for any $x\in X$.
\end{lem}

\begin{rem}
Assume $Y$ is $1$--stable, i.e. $y=S(y_{(-1)})y_{(0)}$ for all $y\in Y$.
Now,
\begin{align*}
p_ni_n & (h^0\otimes\cdots\otimes h^{n-1}\otimes y)
 = p_n(h^0\otimes\cdots\otimes h^{n-1}\otimes y_{(-1)}
     \otimes y_{(0)})\\
 = & \left(S(y_{(-1)(n+1)})y_{(0)(-n)}h^0\otimes\cdots\otimes
           S(y_{(-1)(2)})y_{(0)(-1)}h^{n-1}\otimes
           S(y_{(-1)(1)})y_{(0)(0)}\right)\\
 = & (h^1\otimes\cdots\otimes h^{n-1}\otimes y)
\end{align*}
for all $(h^0\otimes\cdots\otimes h^{n-1}\otimes y)$ from $\B{CM}_n(H,Y)$
which implies $p_*$ is an epimorphism of graded $k$--modules.
\end{rem}

\begin{defn}
Now, define an operator $t_n$ on each $\B{CM}_n(H,Y)$ for $n>0$ by
letting
\begin{align*}
t_n^{-1}(h^0\otimes & \cdots\otimes h^{n-1}\otimes y)\\
 := &  p_n\tau_n^{-1}i_n(h^0\otimes\cdots\otimes h^{n-1}\otimes y)\\
  = & p_n(h^1\otimes\cdots\otimes h^{n-1}\otimes y_{(-2)}\otimes
       y_{(-1)}h^0\otimes y_{(0)})\\
  = & \left(S(h^0_{(n+1)})S(y_{(-1)(n+1)})y_{(0)(-n)}h^1\otimes\cdots
           \otimes S(h^0_{(3)})S(y_{(-1)(3)})y_{(0)(-2)}h^{n-1}\otimes\right.\\
    &  \hspace{1cm}\left.S(h^0_{(2)})S(y_{(-1)(2)})y_{(0)(-1)}y_{(-2)}\otimes
           S(h^0_{(1)})S(y_{(-1)(1)})y_{(0)(0)}\right)\\
  = & \left(S(h^0_{(n+1)})h^1\otimes\cdots
           \otimes S(h^0_{(3)})h^{n-1}\otimes
           S(h^0_{(2)})y_{(-2)}\otimes
           S(h^0_{(1)})S(y_{(-1)})y_{(0)}\right)
\end{align*}
By using the stability again, we get
\begin{align}\label{cm_cyclic}
t_n^{-1}(h^0\otimes\cdots\otimes h^{n-1}\otimes y)
 = \left(S(h^0_{(n+1)})h^1\otimes\cdots
           \otimes S(h^0_{(3)})h^{n-1}\otimes S(h^0_{(2)})y_{(-1)}\otimes
           S(h^0_{(1)})y_{(0)}\right)
\end{align}
\end{defn}

\begin{thm}\label{epimorphism}
Let $H$ be a Hopf algebra and let $Y$ be a stable anti-Yetter-Drinfeld
module.  Then there is an epimorphism of para-cocyclic modules of the
form $\B{T}_*(H,Y)\xra{p_*}\B{CM}_*(H,Y)$.
\end{thm}
\begin{proof}
With these $t_*$ operators at hand for $n=1$, we have
\begin{align*}
t_1^{-1}p_1\tau_1 & (h^0\otimes h^1\otimes y)
 = t_1^{-1}p_1\left(S^{-1}(y_{(-1)})h^1\otimes h^0\otimes y_{(0)}
     \right)\\
 = & t_1^{-1}\left(S(h^0_{(2)})y_{(0)(-1)}S^{-1}(y_{(-1)})h^1
                      \otimes S(h^0_{(1)})y_{(0)(0)}\right)\\
 = & t_1^{-1}\left(S(h^0_{(2)})h^1\otimes S(h^0_{(1)})y\right)\\
 = & \left(S(h^1_{(2)})S^2(h^0_{(2)(1)})S(h^0_{(1)(3)})y_{(-1)}
           h^0_{(1)(1)}\otimes S(h^1_{(1)})S^2(h^0_{(2)(2)})
           S(h^0_{(1)(2)})y_{(0)}\right)\\
 = & \left(S(h^1_{(2)})y_{(-1)}h^0\otimes S(h^1_{(1)})y_{(0)}\right)\\
 = & p_1(h^0\otimes h^1\otimes y)
\end{align*}
On the other hand for $n>1$, we have
\begin{align*}
t_n^{-1}p_n\tau_n & (h^0\otimes\cdots\otimes h^n\otimes y)
 = t_n^{-1}p_n(S^{-1}(y_{(-1)})h^n\otimes h^0\otimes
                    \cdots\otimes h^{n-1}\otimes y_{(0)})\\
 = & t_n^{-1}\left(S(h^{n-1}_{(n+1)})y_{(0)(-n)}S^{-1}(y_{(-1)})h^n
      \otimes S(h^{n-1}_{(n)})y_{(0)(-n+1)}h^0
      \otimes\cdots\otimes S(h^{n-1}_{(1)})y_{(0)(0)}\right)\\
 = & t_n^{-1}\left(S(h^{n-1}_{(n+1)})h^n
      \otimes S(h^{n-1}_{(n)})y_{(-n)}h^0
     \otimes\cdots\otimes S(h^{n-1}_{(1)})y_{(0)}\right)\\
 = & \left(S(h^n_{(n+1)})S^2(h^{n-1}_{(n+1)(1)})
           S(h^{n-1}_{(n)})y_{(-n)}h^0\otimes\cdots\right.\\
   &       \otimes S(h^n_{(2)})S^2(h^{n-1}_{(n+1)(n)})
           S(h^{n-1}_{(1)(3)})y_{(0)(-1)}h^{n-1}_{(1)(1)}
           \otimes S(h^n_{(1)})S^2(h^{n-1}_{(n+1)(n+1)})
           S(h^{n-1}_{(1)(2)})y_{(0)(0)}\\
 = & \left(S(h^n_{(n+1)})y_{(-n)}h^0\otimes\cdots\otimes
           S(h^n_{(2)})y_{(-1)}h^{n-1}\otimes 
           S(h^n_{(1)})y_{(0)}\right)\\
 = & p_n(h^1\otimes\cdots\otimes h^n\otimes y)
\end{align*}
which implies
\begin{align*}
t_np_n=p_n\tau_n
\end{align*}
for all $n\geq 1$.  Notice that
\begin{align*}
p_n\tau_n^{-n-1} & (h^0\otimes\cdots\otimes h^n\otimes y)
 =  p_n\left(y_{(-n-1)})h^0\otimes\cdots\otimes 
              y_{(-1)}h^n\otimes y_{(0)}\right)\\
 = & \left(S(h^n_{(n+1)})S(y_{(-1)(n+1)})y_{(0)(-n)}y_{(-n-1)}h^0\otimes
           \cdots\right.\\
   &       \left.\otimes
           S(h^n_{(2)})S(y_{(-1)(2)})y_{(0)(-1)}y_{(-2)}h^{n-1}
           \otimes S(h^n_{(1)})S(y_{(-1)(1)})y_{(0)(0)}
           \right)\\
 = & \left(S(h^n_{(n+1)})y_{(-n-1)}h^0\otimes\cdots
           \otimes S(h^n_{(2)})y_{(-2)}h^{n-1}
           \otimes S(h^n_{(1)})S(y_{(-1)})y_{(0)}
           \right)
\end{align*}
By using stability, we get
\begin{align*}
 = & \left(S(h^n_{(n+1)})y_{(-n)}h^1\otimes\cdots
           \otimes S(h^n_{(2)})y_{(-1)}h^{n-1}
           \otimes S(h^n_{(1)})y_{(0)}
           \right)\\
 = & p_n(h^1\otimes\cdots\otimes h^n\otimes y)\\
 = & t_n^{-n-1}p_n(h^1\otimes\cdots\otimes h^n\otimes y)
\end{align*}
which means the action of $t_n$ on $\B{CM}_n(H,Y)$ is cyclic of order
$n+1$.

Now let 
\begin{align*}
d_0(h^0\otimes\cdots\otimes h^{n-1}\otimes y)
 = \begin{cases}
   (y_{(-1)}\otimes y_{(0)})
		& \text{ if } n=0\\
   (h^0_{(1)}\otimes h^0_{(2)}\otimes\cdots\otimes y)
		& \text{ if } n>0
   \end{cases}
\end{align*}
and for $0\leq j\leq n+1$ let
\begin{align*}
d_j = & t_{n+1}^jd_0t_n^{-j}
\end{align*}
We are going to show that $d_jp_n=p_{n+1}\partial_j$ for
all $0\leq j\leq n+1$.  Consider the case $j=0$ where
\begin{align*}
p_{n+1}\partial_0 & (h^0\otimes\cdots\otimes h^n\otimes y)\\
 = & p_{n+1}(h^0_{(1)}\otimes h^0_{(2)}\otimes\cdots\otimes y)\\
 = & \begin{cases}
     \left(S(h^0_{(2)(2)})y_{(-1)}h^0_{(1)}\otimes 
           S(h^0_{(2)(1)})y_{(0)}\right)
		& \text{ if } n=0\\
     \left(S(h^n_{(n+1)})y_{(-n)}h^0_{(1)}\otimes
           S(h^n_{(n)})y_{(-n+1)}h^0_{(2)}\otimes\cdots\otimes
           S(h^n_{(1)})y_{(0)}\right)
		& \text{ if } n>0
     \end{cases}\\
 = & \begin{cases}
     \left(S(h^0_{(3)})y_{(-1)}h^0_{(1)}\otimes 
           S(h^0_{(2)})y_{(0)}\right)
		& \text{ if } n=0\\
     \left(S(h^n_{(n)(2)})y_{(-n+1)(1)}h^0_{(1)}\otimes
           S(h^n_{(n)(1)})y_{(-n+1)(2)}h^0_{(2)}\otimes\cdots\otimes
           S(h^n_{(1)})y_{(0)}\right)
		& \text{ if } n>0
     \end{cases}\\
 = & d_0p_n(h^0\otimes\cdots\otimes h^n\otimes y)
\end{align*}
happens only if $Y$ is an anti-Yetter-Drinfeld module.  Then for all
$0\leq j\leq n+1$
\begin{align*}
p_{n+1}\partial_j
 = & p_{n+1}\tau_{n+1}^j\partial_0\tau_n^{-j}
 =  t_{n+1}^jp_{n+1}\partial_0\tau_n^{-j}
 =  t_{n+1}^jd_0p_n\tau_n^{-j}
 =  t_{n+1}^jd_0t_n^{-j}p_n
 =  d_jp_n
\end{align*}
as we wanted to show.
\end{proof}

\begin{thm}\label{shift}
Let $H$ be Hopf algebra and let $Y$ be a stable anti-Yetter-Drinfeld
module.  Then $\B{CM}_*(H,Y)$ is isomorphic to the Connes--Moscovici
cyclic complex.
\end{thm}
\begin{proof}
The fact that $\B{CM}_*(H,Y)$ is a cocyclic module follows from
Theorem~\ref{epimorphism}.  Now, let us see the cosimplicial structure
maps on $\B{CM}_*(H,Y)$ explicitly: Take $(h^0\otimes\cdots\otimes
h^{n-1}\otimes y)$ from $\B{CM}_n(H,Y)$.  Since
$d_jp_n=p_{n+1}\partial_j$ for any $0\leq j\leq n+1$ and $p_ni_n=id_n$ I
have $d_j=p_{n+1}\partial_ji_n$.  Therefore for $n=0$,
\begin{align*}
d_j(y)
 = & \begin{cases}
     p_{n+1}(y_{(-2)}\otimes y_{(-1)}\otimes y_{(0)})
        & \text{ if } j=0\\
     p_{n+1}(y_{(-2)}\otimes y_{(-1)}y_{(-3)}\otimes y_{(0)})
        & \text{ if } j=1
     \end{cases}\\
 = & \begin{cases}
     (S(y_{(-2)})y_{(-1)}y_{(-4)}\otimes S(y_{(-3)})y_{(0)})
        & \text{ if } j=0\\
     (S(y_{(-5)})S(y_{(-2)})y_{(-1)}y_{(-4)}
        \otimes S(y_{(-6)})S(y_{(-3)})y_{(0)})
        & \text{ if } j=1
     \end{cases}\\
  = & \begin{cases}
      (y_{(-1)}\otimes y_{(0)})
        & \text{ if } j=0\\
      (1\otimes y) 
        & \text{ if } j=1
      \end{cases}
\end{align*}
For $n>0$,
\begin{align*}
d_j(h^0\otimes\cdots\otimes h^{n-1}\otimes y)
 = & \begin{cases}
     p_{n+1}(\cdots\otimes h^j_{(1)}\otimes
             h^j_{(2)}\otimes\cdots\otimes y_{(-1)}\otimes y_{(0)})
        & \text{ if } 0\leq j\leq n-1\\
     p_{n+1}(h^0\otimes\cdots\otimes h^{n-1}\otimes y_{(-2)}
             \otimes y_{(-1)}\otimes y_{(0)})
        & \text{ if } j=n\\
     p_{n+1}(h^0_{(2)}\otimes h^1\otimes\cdots\otimes h^{n-1}
             \otimes y_{(-2)}\otimes y_{(-1)}h^0_{(1)}\otimes y_{(0)})
        & \text{ if } j=n+1
     \end{cases}
\end{align*}
Let me investigate these cases separately: For $0\leq j\leq n-1$ one has
\begin{align*}
p_{n+1} & (\cdots\otimes h^j_{(1)}\otimes
             h^j_{(2)}\otimes\cdots\otimes y_{(-1)}\otimes y_{(0)})\\
 = & \left(S(y_{(-1)(n+2)})y_{(0)(-n-1)}h^0\otimes\cdots\otimes
           S(y_{(-1)(n+2-j)})y_{(0)(-n-1+j)}h^j_{(1)}\otimes\right.\\
   & \hspace{1.5cm}\left.
           S(y_{(-1)(n+1-j)})y_{(0)(-n+j)}h^j_{(2)}\otimes\cdots\otimes
           S(y_{(-1)(2)})y_{(0)(-1)}h^{n-1}\otimes
           S(y_{(-1)(1)})y_{(0)(0)}\right)\\
 = & (\cdots\otimes h^j_{(1)}\otimes h^j_{(2)}\otimes\cdots\otimes y)  
\end{align*}
For $j=n$ one has
\begin{align*}
p_{n+1} & (h^0\otimes\cdots\otimes h^{n-1}\otimes y_{(-2)}
             \otimes y_{(-1)}\otimes y_{(0)})\\
 = & \left(S(y_{(-1)(n+2)})y_{(0)(-n-1)}h^0\otimes\cdots\otimes 
           S(y_{(-1)(3)})y_{(0)(-2)}h^{n-1}\otimes\right.\\
   & \hspace{1.5cm}\left.
           S(y_{(-1)(2)})y_{(0)(-1)}y_{(-2)}\otimes 
           S(y_{(-1)(1)})y_{(0)(0)}\right)\\
 = & (h^0\otimes\cdots\otimes h^{n-1}\otimes y_{(-1)}\otimes y_{(0)})
\end{align*}
And finally for $j=n+1$ one has
\begin{align*}
p_{n+1} & (h^0_{(2)}\otimes h^1\otimes\cdots\otimes h^{n-1}
        \otimes y_{(-2)}\otimes y_{(-1)}h^0_{(1)}\otimes y_{(0)})\\
 = & \left(1\otimes S(h^0_{(n+1)})h^1\otimes\cdots\otimes
            S(h^0_{(3)})h^{n-1}\otimes S(h^0_{(2)})y_{(-1)}
            \otimes S(h^0_{(1)})y_{(0)}\right)
\end{align*}
which are slightly different than the face maps defined in
\cite{ConnesMoscovici:HopfCyclicCohomology} and
\cite{Crainic:CyclicCohomologyOfHopfAlgebras}.  The cosimplicial
structure morphisms on the Connes--Moscovici cyclic complex
$\C{C}_*(H,Y)$ are
\begin{align*}
\widetilde{d}_j(h^0\otimes\cdots\otimes h^{n-1}\otimes y)
 = & \begin{cases}
     (1\otimes h^0\otimes\cdots\otimes h^{n-1}\otimes y)
        & \text{ if } j=0\\
     (\cdots\otimes h^{j-1}_{(1)}\otimes h^{j-1}_{(2)}\otimes
      \cdots\otimes y)
        & \text{ if } 1\leq j\leq n\\
     (h^0\otimes\cdots\otimes h^{n-1}\otimes y_{(-1)}\otimes y_{(0)})
        & \text{ if } j=n+1
     \end{cases}
\end{align*}
We define the isomorphism $\B{CM}_*(H,Y)\xra{\alpha_*}\C{C}_*(H,Y)$ as
\begin{align*}
\alpha_n({\bf h}\otimes y) = t_n^{-1}({\bf h}\otimes y)
\end{align*}
for all ${\bf h}\otimes y$ from $\B{CM}_*(H,Y)$.  Obviously
\begin{align*}
\alpha_nt_n = & t_nt_n^{-1} = t_n^{-1}t_n = t_n\alpha_n
\end{align*}
Note that for $1\leq j\leq n+1$
\begin{align*}
\alpha_{n+1}d_j
 = & t_{n+1}^{-1}d_j
 =   d_{j-1}t_n^{-1}
 =   \widetilde{d}_j\alpha_n
\end{align*}
This leaves
\begin{align*}
\alpha_{n+1}d_0(h^0\otimes\cdots\otimes h^{n-1}\otimes y)
 = & t_{n+1}^{-1}d_0(h^0\otimes\cdots\otimes h^{n-1}\otimes y)\\
 = & d_{n+1}(h^0\otimes\cdots\otimes h^{n-1}\otimes y)\\
 = & \left(1\otimes S(h^0_{(n+1)})h^1\otimes\cdots\otimes
            S(h^0_{(3)})h^{n-1}\otimes S(h^0_{(2)})y_{(-1)}
            \otimes S(h^0_{(1)})y_{(0)}\right)\\
 = & \widetilde{d_0}t_n^{-1}(h^0\otimes\cdots\otimes h^{n-1}\otimes y)\\
 = & \widetilde{d_0}\alpha_n(h^0\otimes\cdots\otimes h^{n-1}\otimes y)
\end{align*}
for all $(h^0\otimes\cdots\otimes h^{n-1}\otimes y)$ from
$\B{CM}_n(H,Y)$ as we wanted to show.
\end{proof}

\begin{thm}\label{factoring_isomorphism}
Let $H$ and $Y$ be as before.  Then $\B{T}_*(H,Y)\xra{p_*}\B{CM}_*(H,Y)$
factors as
\begin{align}\label{factoring}
\B{T}_*(H,Y)\xra{q_*}{}_H\B{T}_*(H,Y)\xra{p'_*}\B{CM}_*(H,Y)
\end{align}
Moreover, $p'_*$ is an isomorphism of cocyclic $k$--modules.
\end{thm}
\begin{proof}
Take $(h^0\otimes\cdots\otimes h^n\otimes y)$ from $\B{T}_*(H,Y)$ and
$h\in H$ and consider
\begin{align*}
p_nL_h & (h^0\otimes\cdots\otimes h^n\otimes y)\\
 = & p_n(h_{(1)}h^0\otimes\cdots\otimes h_{(n+1)}h^n\otimes
         h_{(n+2)}y)\\
 = & \left(S(h^n_{(n+1)})S(h_{(n+1)(n+1)})
           h_{(n+2)(1)}y_{(-n)}S^{-1}(h_{(n+2)(2n+1)})h_{(1)}h^0
           \otimes\cdots\otimes\right.\\
   & \hspace{1cm}
           S(h^n_{(2)})S(h_{(n+1)(2)})
           h_{(n+2)(n)}y_{(-1)}S^{-1}(h_{(n+2)(n+2)})h_{(n)}h^{n-1}\otimes\\
   & \hspace{1cm}\left. S(h^n_{(1)})S(h_{(n+1)(1)})
           h_{(n+2)(n+1)}y_{(0)} \right)\\
 = & \left(S(h^n_{(n+1)})y_{(-n)}S^{-1}(h_{(n+2)(n+1)})h_{(1)}h^0
           \otimes\cdots\otimes\right.\\
   & \hspace{1cm}\left.
           S(h^n_{(2)})y_{(-1)}S^{-1}(h_{(n+2)(2)})h_{(n)}h^{n-1}
           \otimes S(h^n_{(1)})S(h_{(n+1)})
           h_{(n+2)(1)}y_{(0)} \right)\\
 = & \epsilon(h)\left(S(h^n_{(n+1)})y_{(-n)}h^0
           \otimes\cdots\otimes S(h^n_{(2)})y_{(-1)}h^{n-1}
           \otimes S(h^n_{(1)})y_{(0)} \right)
\end{align*}
This shows that one has a factoring of the form~(\ref{factoring})

Let $H\otimes\B{T}_*(H,Y)\xra{\rho^{H,Y}}\B{T}_*(H,Y)$ denote the left
$H$ action and let $\epsilon$ denote the trivial action.  Consider the
following commutative diagram
\begin{align*}
\begin{CD}
ker(p_*)            @>>> \B{T}_*(H,Y) @>{p_*}>> \B{CM}_*(H,Y)\\
@A{id}AA            @A{id}AA                   @AA{p'_*}A\\
im(\epsilon-\rho^{H,Y}) @>>> \B{T}_*(H,Y) @>>{q_*}> {}_H\B{T}_*(H,Y)
\end{CD}
\end{align*}
We need to show that any element in $ker(p_*)$ is in
$im(\epsilon-\rho^{H,Y})$.

Since $p_*i_*=id_*$, we have $ker(p_*)=im(id_*-i_*p_*)$.  Then take
$(id_*-i_*p_*)({\bf h}\otimes y)$ from $ker(p_*)$ and consider
\begin{align*}
i_np_n & (h^0\otimes\cdots\otimes h^n\otimes y)
 = i_n\left(S(h^n_{(n+1)})y_{(-n)}h^0\otimes\cdots\otimes
              S(h^n_{(2)})y_{(-1)}h^{n-1}\otimes S(h^n_{(1)})y_{(0)}\right)\\
 = & \left(S(h^n_{(n+1)})y_{(-n)}h^0\otimes\cdots\otimes
              S(h^n_{(2)})y_{(-1)}h^{n-1}\otimes 
              S(h^n_{(1)(3)})y_{(0)(-1)}h^n_{(1)(1)}\otimes
              S(h^n_{(1)(2)})y_{(0)(0)}\right)\\
 = & \left(S(h^n_{(n+3)})y_{(-n-2)}h^0\otimes\cdots\otimes
              S(h^n_{(4)})y_{(-3)}h^{n-1}\otimes 
              S(h^n_{(3)})y_{(-2)}h^n_{(1)}\otimes
              S(h^n_{(2)})y_{(-1)}y_{(0)}\right)\\ 
 = & L_{S(h^n_{(2)})y_{(-1)}}
       (h^0\otimes\cdots\otimes h^{n-1}\otimes h^n_{(1)}\otimes y_{(0)})
\end{align*}
Since $\epsilon(S(h_{(1)}))h_{(2)}=h$ for any $h\in H$, I see that
\begin{align*}
(id_n-i_np_n)(h^0\otimes\cdots\otimes h^n)
 = \left(\epsilon(S(h^n_{(2)})y_{(-1)})-S(h^n_{(2)})y_{(-1)}\right)
     \cdot(h^0\otimes\cdots\otimes h^{n-1}\otimes h^n_{(1)}\otimes y_{(0)})  
\end{align*}
which means $ker(p_*)\subseteq im(\epsilon-\rho^{H,Y})$ as I wanted to show.
\end{proof}

\begin{rem}
I would like to stress that $\B{T}_*(H,Y)$ IS NOT a para-cocyclic
$H$--module since the morphisms $\partial_{n+1}$ and the action of any
$h\in H$ on $\B{T}_n(H,Y)$ may not commute for all $n\geq 1$ unless $H$
is cocommutative.  However, ${}_H\B{T}_*(H,Y)$ is a cosimplicial
$H$--module, since the action of $H$ is trivial.
\end{rem}

\section{A new cyclic homology for bialgebras}\label{NewHomology}

In this section, we assume $Y$ is just a stable $H$--module/comodule.

\begin{defn}
One can see that
\begin{align}\label{null}
p_nL_h\tau_n^i
 = & \epsilon(S^{-1}(h))p_n\tau_n^i
 = \epsilon(S^{-1}(h))t_n^ip_n
 = t_n^ip_n L_h
 = p_n\tau_n^iL_h
\end{align}
For all $0\leq j\leq n+1$, $i\in\B{Z}$, $h\in H$ and $({\bf h}\otimes
y)\in\B{T}_n(H,Y)$.  Now, let $I_{[H,\C{C}]}$ be the $k$--submodule of
$\B{T}_*(H,Y)$ containing expressions of the form
\begin{align*}
[L_h,\tau_n^j]({\bf h}\otimes y)
 = & L_h\tau_n^j({\bf h}\otimes y) - \tau_n^jL_h({\bf h}\otimes y)
\end{align*}
for all $({\bf h}\otimes y)$ from $H^{\otimes y}\otimes Y$, $h$ from
$H$, $j\in\B{Z}$ and $n\geq 0$.
\end{defn}

\begin{thm}\label{main_isomorphism}
$\B{T}_*(H,Y)/I_{[\C{C},H]}$ is a para-cocyclic $H$--module and
$\B{CM}_*(H,Y) \cong {}_H\left(\B{T}_*(H,Y)/I_{[\C{C},H]}\right)$.
\end{thm}
\begin{proof}
$I_{[H,\C{C}]}$ is contained in $ker(p_*)$ since I have (\ref{null}).
Note that $I_{[H,\C{C}]}$ is a graded $H$--submodule since
\begin{align*}
L_x[L_h,\tau_n^i]
 = & L_xL_h\tau_n^i - L_x\tau_n^iL_h
 =   L_{xh}\tau_n^i - \tau_n^iL_{xh} + \tau_n^iL_xL_h - L_x\tau_n^iL_h\\
 = & [L_{xh},\tau_n^i]+ [\tau_n^i,L_x]L_h
\end{align*}
for all $x,h\in H$ and $i\in \B{Z}$.  So I get a sequence of
epimorphism of graded $H$--modules of the form
\begin{align*}
\B{T}_*(H,Y)\xra{q^1_*}\B{T}_*(H,Y)/I_{[H,\C{C}]}\xra{q^2_*}{}_H\B{T}_*(H,Y)
\end{align*}
Since the functor $X\mapsto {}_H X$ is right exact and these morphisms are
graded $H$--module morphisms, I obtain another sequence of epimorphisms
\begin{align*}
{}_H\B{T}_*(H,Y)\xra{{}_H(q^1_*)}
  {}_H\left(\B{T}_*(H,Y)/I_{[H,\C{C}]}\right)
  \xra{{}_H(q^2_*)}{}_H\B{T}_*(H,Y)
\end{align*}
where the composition is identity.  Then ${}_H(q^1_*)$ is both an epi-
and a mono-morphism, i.e.
\begin{align*}
{}_H\left(\B{T}_*(H,Y)/I_{[H,\C{C}]}\right)\xra{{}_H(q^2_*)}
{}_H\B{T}_*(H,Y)
\end{align*}
is an isomorphism.

Notice that $I_{[H,\C{C}]}$ is stable under the action of $\tau_*$ since
\begin{align*}
\tau_n^i[L_h,\tau_n^j]
 = & \tau_n^i L_h \tau_n^j - \tau_n^{i+j}L_h 
 =   \tau_n^i L_h \tau_n^j - L_h\tau_n^{i+j} + L_h\tau_n^{i+j}
    - \tau_n^{i+j}L_h
 = [\tau_n^i,L_h]\tau_n^j + [L_n,\tau_n^{i+j}]
\end{align*}
Furthermore consider,
\begin{align}\label{here}
\partial_0[L_h,\tau_*^{-j}]
 = & -[\partial_0,\tau_*^{-j}]L_h + [\partial_0 L_h,\tau_*^{-j}]
\end{align}
by defining
\begin{align*}
[\partial_0,\tau_*^{-j}]
  :=\partial_0\tau_n^{-j}  - \tau_{n+1}^{-j}\partial_0
   =\tau_{n+1}^{-j}\left(\partial_j - \partial_0\right)
\end{align*}
Notice also that $[L_x,\partial_0]=0$ since
\begin{align}\label{partial}
L_x\partial_0(h^0\otimes\cdots\otimes h^n\otimes y)
 = & \left(x_{(1)}h^0_{(1)}\otimes x_{(2)}h^0_{(2)}\otimes
    \cdots\otimes x_{(n+2)}h^n\otimes x_{(n+3)}y\right)\\
 = & \partial_0L_x(h^0\otimes\cdots\otimes h^n\otimes y)
\end{align}
Then (\ref{here}) reads as
\begin{align*}
\partial_0[L_h,\tau_*^{-j}]
 = & -\tau_*^{-j}(\partial_j-\partial_0)L_h 
     + [L_h\partial_0 ,\tau_*^{-j}]\\
 = & -\tau_*^{-j}(\partial_j-\partial_0)L_h 
     + [L_h,\tau_*^{-j}]\partial_0 + L_h[\partial_0,\tau_*^{-j}]\\
 = & [L_h,\tau_*^{-j}(\partial_j-\partial_0)]
     + [L_h,\tau_*^{-j}]\partial_0\\
 = & [L_h,\tau_*^{-j}](\partial_j-\partial_0)
     + \tau_*^{-j}[L_h,(\partial_j-\partial_0)]
     + [L_h,\tau_*^{-j}]\partial_0\\
 = & [L_h,\tau_*^{-j}]\partial_j
     + \tau_*^{-j}[L_h,\partial_j]
\end{align*}
So, the problem of showing $I_{[H,\C{C}]}$ is stable under the action of
$\partial_0$ reduces to showing $[L_h,\partial_j]({\bf h}\otimes y)\in
I_{[H,\C{C}]}$ for all $h\in H$, $({\bf h}\otimes
y)\in\B{T}_n(H,Y)$ and for all $0\leq j\leq n+1$.  I have shown
above in (\ref{partial}) that $[L_h,\partial_0]({\bf h}\otimes y)=0$ for
all $({\bf h}\otimes y)\in \B{T}_n(H,Y)$.  Similarly, one can show for
$0\leq j\leq n$ that $[L_h,\partial_j]({\bf h}\otimes y)=0$.  That
leaves out
\begin{align*}
[L_h,\partial_{n+1}]({\bf h}\otimes y)
 = & [L_h,\tau_{n+1}^{-1}\partial_0]({\bf h}\otimes y)\\
 = & [L_h,\tau_{n+1}^{-1}]\partial_0({\bf h}\otimes y)
     + \tau_{n+1}[L_h,\partial_0]({\bf h},y)\\
 = & [L_h,\tau_{n+1}^{-1}]\partial_0({\bf h}\otimes y)
\end{align*}
which is inside $I_{[H,\C{C}]}$.  So, I showed that $I_{[H,\C{C}]}$ is
stable under the actions of $\tau_*$ and $\partial_0$ which implies
$I_{[H,\C{C}]}$ is a para-cocyclic submodule of $\B{T}_*(H,Y)$, which
in turn implies $\B{T}_*(H,Y)/I_{[H,\C{C}]}$ is a para-cocyclic
module.  Moreover, $\B{T}_*(H,Y)/I_{[H,\C{C}]}$ is a para-cocyclic
$H$--module since $[L_h,\partial_j]=0$ for all $h\in H$ on
$\B{T}_*(H,Y)/I_{[H,\C{C}]}$.
\end{proof}

\begin{rem}
The most important feature of the para-cocyclic $H$--module
$\B{T}_*(H,Y)/I_{[H,\C{C}]}$ is that it can be defined without requiring
$H$ to be a Hopf algebra or $Y$ to be an anti-Yetter-Drinfeld
$H$--module.  And Theorem~\ref{main_isomorphism} tells us that in the
case of $H$ is a Hopf algebra and $Y$ is an anti-Yetter-Drinfeld
$H$--module, the Connes--Moscovici cyclic homology of $H$ with
coefficients in $Y$ can be recovered as the homology of the
$H$--coinvariants of the complex $\B{T}_*(H,Y)/I_{[H,\C{C}]}$.
\end{rem}

Assume $B$ is an ordinary unital/counital associative/coassociative
bialgebra.  Let $X$ be a (left) $B$--module coalgebra and $Y$ be a
left $B$--module/comodule which is $0$--stable, i.e. 
\begin{align*}
\sum_y y_{(-1)}y_{(0)} = y
\end{align*}
for all $y\in Y$.  Let $\B{T}_*(X,Y)=\left\{X^{\otimes n+1}\otimes
Y\right\}_{n\geq 0}$ be the graded $B$--module with diagonal action
\begin{align*}
L_b(x^0\otimes\cdots\otimes x^n\otimes y)
 := \left(b_{(1)} x^0\otimes\cdots\otimes b_{(n+1)}x^n\otimes
          b_{(n+1)}y\right)
\end{align*}
for all $b\in B$ and $({\bf x}\otimes y)=(x^0\otimes\cdots\otimes
x^n\otimes y)$ in $\B{T}_n(X,Y)$ and $n\geq 0$.  Define a cosimplicial
structure by
\begin{align*}
\partial_j(x^0\otimes\cdots\otimes x^n\otimes y)
 = \begin{cases}
   (\cdots\otimes x^j_{(1)}\otimes x^j_{(2)}\otimes\cdots\otimes y)
		& \text{ if } 0\leq j\leq n\\
   (x^0_{(2)}\otimes x^1\otimes\cdots\otimes y_{(-1)}x^0_{(1)}
    \otimes y_{(0)})
		& \text{ if } j=n+1
   \end{cases}
\end{align*}
Define also
\begin{align*}
\tau_n^{-1}(x^0\otimes\cdots\otimes x^n\otimes y)
 = (x^1\otimes\cdots\otimes x^n\otimes y_{(-1)}x^0\otimes y_{(0)})
\end{align*}
which is not necessarily invertible.  If $B$ is a Hopf algebra, then
$\tau_n^{-1}$ certainly is invertible and together with the cosimplicial
maps defines a para-cocyclic structure.  Let $I_{[B,\C{C}]}$ be the
graded $k$--submodule of $\B{T}_*(X,Y)$ which is generated by elements
of the form
\begin{align*}
[L_b,\tau_n^{-j}]({\bf x}\otimes y)
\end{align*}
for $b\in B$, $({\bf x}\otimes y)\in \B{T}_n(X,Y)$ and $n\geq 0$.  As
before, $I_{[B,\C{C}]}$ is a graded $B$--submodule of $\B{T}_*(X,Y)$,
and
\begin{align*}
[\partial_0,\tau_*^{-j}]
 =\partial_0\tau_n^{-j}-\tau_{n+1}^{-j}\partial_0
 = \tau_{n+1}^{-j}(\partial_j-\partial_0)
\end{align*}
Therefore $I_{[B,\C{C}]}$ is a cosimplicial submodule of $\B{T}_*(X,Y)$
which is stable under the actions of $\tau_*^{-1}$.  This makes
$\B{T}_*(X,Y)/I_{[B,\C{C}]}$ a cosimplicial $B$--module which is short
of being a para-cocyclic $B$--module since $\tau_*^{-1}$ may not be
invertible.  However, if $B$ is a Hopf algebra, then $\tau_*^{-1}$ is
invertible.  Moreover, even in the case of $B$ is an ordinary bialgebra,
the quotient ${}_B\left(\B{T}_*(X,Y)/I_{[B,\C{C}]}\right)$ is still a
cocyclic module since
\begin{align*}
t_n^{-n-1}[x^0\otimes\cdots\otimes x^n\otimes y]
 = & [y_{(-n-1)}x^0\otimes\cdots\otimes y_{(-1)}x^n\otimes y_{(0)}]\\
 = & [y_{(-n-2)}x^0\otimes\cdots\otimes y_{(-2)}x^n\otimes y_{(-1)}y_{(0)}]\\
 = & y_{(-1)}[x^0\otimes\cdots\otimes x^n\otimes y_{(0)}]\\
 = & \epsilon(y_{(-1)})[x^0\otimes\cdots\otimes x^n\otimes y_{(0)}]\\
 = & [x^0\otimes\cdots\otimes x^n\otimes y]
\end{align*}

For simplicity, denote the para-cocyclic $B$--module
$\B{T}_*(X,Y)/I_{[B,\C{C}]}$ by $\B{PCM}_*(X,Y)$.  So, the result I
have above reads as follows:

\begin{thm}\label{generalization}
Let $B$ be an ordinary associative/coassociative unital/counital
bialgebra.  Let $X$ be a $B$--module coalgebra and let $Y$ be a
$0$--stable $B$--module/comodule.  Then the graded $B$--module
$\B{PCM}_*(X,Y)$ defined above is a cosimplicial $B$--module which is
short of being a para-cocyclic module since $\tau_*^{-1}$ may not be
invertible.  However, if $B$ is a Hopf algebra then $\B{PCM}_*(X,Y)$ is
para-cocyclic $B$--module.  Moreover, regardless of $B$ being a Hopf
algebra, $\B{CM}_*(X,Y):={}_B\B{PCM}_*(X,Y)$ is {\bf always} a cocyclic
module.  Finally, if $X=B$ is a Hopf algebra and $Y$ is a stable
anti-Yetter-Drinfeld module then cyclic homology of the cyclic complex
$\B{CM}_*(X,Y)$ computes the Connes--Moscovici cyclic homology of $B$
with coefficients in $Y$.  We denote this new homology by
$HC^\B{CM}_*(X,Y)$.
\end{thm}

\begin{rem}
The functor $HC^\B{CM}_*(\cdot,\cdot)$ which computes the bialgebra
cyclic homology of a pair $(X,Y)$ where $X$ is a $B$--module coalgebra
and $Y$ is a $0$--stable $B$--module/comodule is covariant in both
variables.  Since we didn't introduced the codegeneracy maps, the reader
should take caution in computing the cyclic homology of the cyclic
complex $\B{CM}_*(X,Y)$.  One can use $(b,B)$ complex (2.1.7 to 2.1.11
in \cite{Loday:CyclicHomology}) if one proves the results obtained so
far with the codegeneracy morphisms.  If one assumes $k$ is of
characteristic $0$, then one can use cyclic invariants of the Hochschild
complex of $\B{CM}_*(X,Y)$, as it is defined, to compute the cyclic
homology (2.1.4 and 2.1.5 in \cite{Loday:CyclicHomology}.)  If the
reader wants a field with arbitrary characteristic, he/she can compute
the cyclic homology of $\B{CM}_*(X,Y)$ by using the cyclic bicomplex
(2.1.2 in \cite{Loday:CyclicHomology}.)
\end{rem}

\begin{cor}
Assume $H$ is a cocommutative Hopf algebra and $Y$ is a stable
anti-Yetter-Drinfeld module.  Then $\B{PCM}_*(H,Y)=\B{T}_*(H,Y)$.
\end{cor}
\begin{proof}
$\B{PCM}_*(H,Y)$ is obtained from $\B{T}_*(H,Y)$ by dividing out the
para-cocyclic submodule generated by elements of the form
\begin{align*}
[\tau_n^j,L_x] & (h^0\otimes\cdots\otimes h^n\otimes y)\\
 = & \tau_n^j\left(x_{(1)}h^0\otimes\cdots\otimes x_{(n+1)}h^n\otimes
                 x_{(n+2)}\cdot y\right)\\
   & - L_x\left(S^{-1}(y_{(-1)})h^{n+1-j}\otimes\cdots\otimes
                S^{-1}(y_{(-j)})h^n\otimes h^0\otimes\cdots\otimes
                h^{n-j}\otimes y_{(0)}\right)\\
 = & \left(x_{(n+2)}y_{(-1)}S^{-1}(x_{(n+2j+2)})x_{(n+2-j)}h^{n+1-j}
                \otimes\cdots\otimes
                x_{(n+j+1)}y_{(-j)}S^{-1}(x_{(n+j+3)})x_{(n+1)}
                (y_{(-j)})h^n\right.\\
   & \hspace{1.5cm}\left.
                \otimes x_{(1)}h^0\otimes\cdots\otimes
                x_{(n+1-j)}h^{n-j}\otimes x_{(n+j+2)}y_{(0)}\right)\\
   & - L_x\left(S^{-1}(y_{(-1)})h^{n+1-j}\otimes\cdots\otimes
                S^{-1}(y_{(-j)})h^n\otimes h^0\otimes\cdots\otimes
                h^{n-j}\otimes y_{(0)}\right)\\
 = & 0
\end{align*}
since $H$ is cocommutative, for any $j\geq 0$.  The proof for $j< 0$ is
similar.
\end{proof}

This means, in the classical cases where $H$ is the group algebra of a
group $G$ or $H$ is the universal enveloping algebra of a Lie algebra
$\G{g}$, if one chooses $Y$ to be stable anti-Yetter-Drinfeld module,
$\B{PCM}_*(H,Y)$ is going to be the same as $\B{T}_*(H,Y)$.

\begin{cor}
Let $H$ and $Y$ be as before.  Then
\begin{align*}
\B{PCM}_*(H,Y)\cong H\otimes \B{CM}_*(H,Y)
\end{align*}
\end{cor}
\begin{proof}
$coad(H)$ is a trivial $H$--comodule since $H$ is a cocommutative
coalgebra. Then
\begin{align*}
\B{PCM}_*(H,Y) = \B{T}_*(H,Y) \cong BC_*(coad(H),H,Y)
 \cong H\otimes BC_*(k,H,Y)
\end{align*}
Since $k$ is a field, $H$ is a projective $k$--module.  Result follows.
\end{proof}

\section{Some useful technical lemmata}\label{TechnicalResults}

In this section, we assume $H$ is a Hopf algebra with an invertible
antipode and $B$ is a unital/counital bialgebra.

\begin{lem}\label{equivariant}
Let $X$ and $Y$ be two left $H$--modules.  Then $X\otimes Y$ is a left
$H$--module via the diagonal action.  Moreover, $X$ can be considered as
a right $H$--module $X^R$ via $x\cdot h:=S^{-1}(h)\cdot x$ for any $x\in
X$ and $h\in H$.  Then ${}_H(X\otimes Y)\cong X^R\eotimes{H}Y$.
\end{lem}

\begin{proof}
The fact that $X\otimes Y$ is a left $H$--module is clear.  Consider the
quotient morphism $(X\otimes Y)\xra{q}(X^R \eotimes{H} Y)$ and I get
\begin{align*}
q\left(h\cdot (x\otimes y)\right)
 = & q\left(h_{(1)}\cdot x\otimes h_{(2)}\cdot y\right)\\
 = & q\left(x\cdot S(h_{(1)})\otimes h_{(2)}\cdot y\right)\\
 = & \epsilon(h)(x\otimes y)
\end{align*}
This means $q$ factors as $(X\otimes Y)\xra{q'}{}_H(X\otimes Y)
\xra{q''}(X^R\eotimes{H}Y)$.  Moreover, if I define $(X^R\eotimes{H}Y)
\xra{u}{}_H(X\otimes Y)$ by letting $u(x\eotimes{H} y)=(x\otimes y)$, then
one can see that
\begin{align*}
u(x\cdot S(h)\otimes y) - u(x\otimes S(h)\cdot y)
 = & (x\cdot S(h)\otimes y) - (x\otimes S(h)\cdot y)\\
 = & (h\cdot x\otimes y) - (x\otimes S(h)\cdot y)\\
 = & \epsilon(S(h_{(1)}))(h_{(2)}\cdot x\otimes y) 
     - (x\otimes S(h)\cdot y)\\
 = & (S(h_{(1)(2)})h_{(2)}\cdot x\otimes S(h_{(1)(1)})\cdot y) 
     - (x\otimes S(h)\cdot y)\\
 = & 0
\end{align*}
which means $u$ is well-defined and is the inverse of $q''$.  This
proves the assertion.
\end{proof}

\begin{lem}\label{quotient}
Let $Y$ be a $0$--stable $B$--module/comodule.  Let $J$ be a coideal and
a left ideal in $B$ such that $p_*(B^{\otimes j}\otimes J\otimes
B^{\otimes n-j})\equiv 0$ for any $0\leq j\leq n$ where
$\B{T}_*(B,Y)\xra{p_*}\B{CM}_*(B,Y)$ is defined earlier.  Then $B/J$ is
a $B$--module coalgebra and $\B{CM}_*(B,Y)\cong {}_B\B{T}_*(B/J,Y)$.
\end{lem}
\begin{proof}
Define
\begin{align*}
\B{T}_*(J,Y)=\left\{\bigoplus_{j=0}^n B^{\otimes j}\otimes J\otimes 
                   B^{n-j}\otimes Y\right\}_{n\geq 0}
\end{align*}
Since $J$ is a coideal, $B/J$ is a coalgebra.  It is clear that $B/J$ is
a left $B$--module and since $B$ itself is a $B$--module coalgebra, the
quotient $B/J$ is naturally a $B$--module coalgebra as well.  Moreover,
since $J$ is a coideal and a left ideal, $\B{T}_*(B,Y)$ is a
para-cocyclic submodule of $\B{T}_*(B,Y)$.  Furthermore, the hypothesis
on $J$ implies $p_*\B{T}_*(J,Y)\equiv 0$ which means $p_*$ factors as
\begin{align*}
\B{T}_*(B,Y)\xra{\pi_*}\B{T}_*(B,Y)/\B{T}_*(J,Y)\xra{\pi'_*}\B{CM}_*(B,Y)
\end{align*}
where the middle term is $\B{T}_*(B/J,Y)$.  Since $p_*$ is an
epimorphism, so is $\pi'_*$.  Then by taking coinvariants with respect
to $B$, I see that ${}_B\B{T}_*(B/J,Y)\cong \B{CM}_*(B,Y)$.  This
finishes the proof.
\end{proof}

\begin{lem}\label{stability1}
Assume $Y$ is an arbitrary $B$--comodule considered as a $B$--module via
the counit.  In other words, let $b\cdot y=\epsilon(b)y$ for any $b\in
B$ and $y\in Y$.  Then $Y$ is a $0$--stable $B$--module/comodule
\end{lem}

\begin{proof}
Clearly, for any $y\in Y$ one has $y_{(-1)}\cdot y_{(0)} =
\epsilon(y_{(-1)})y_{(0)} = y$ by definition of the action.  This means
$Y$ is now a $0$--stable $B$--module/comodule.
\end{proof}

\begin{lem}\label{stability2}
Assume $Y$ is an arbitrary $B$--module considered as a $B$--comodule via
the trivial coaction.  In other words let $\rho_Y(y)=1\otimes y$ for any
$y\in Y$.  Then $Y$ is a $0$--stable $B$--module/comodule.
\end{lem}

\begin{rem}
Although stability of the coefficient module/comodule is enough to
define bialgebra cyclic homology, one might ask a more involved
interplay between the module and comodule structures in order to perform
computations in cyclic homology.  The aYD condition is an avenue one can
follow.  The major classes of stable module/comodules satisfying the aYD
condition are Connes and Moscovici's 1-dimensional module/comodules
coming from a modular pair \cite{ConnesMoscovici:HopfCyclicCohomologyI}
\cite{ConnesMoscovici:CyclicCohomologyOfHopfAlgebras} and Hopf algebras
coming from Hopf--Galois extensions \cite{JaraStefan:HopfGalois}.  More
examples can be found in \cite{Khalkhali:SaYDModules}.  Now,
Lemma~\ref{stability1} and Lemma~\ref{stability2} give us a whole new
class of stable module/comodules which do not satisfy aYD, therefore
suitable only for the extended bialgebra cyclic homology I defined
above.  For instance, any left coideal $Y$ of a bialgebra $B$ with
trivial action can now be considered as a $0$--stable
$B$--module/comodule.  Similarly, any left ideal $Y$ of a bialgebra $B$
with trivial coaction can also be considered as a $0$--stable
$B$--module/comodule.
\end{rem}

\begin{defn}
Define an endomorphism of the graded $k$--module $\B{T}_*(H,Y)$ by
\begin{align*}
\kappa_x &(h^0\otimes\cdots\otimes h^n\otimes y)\nonumber\\
 = & [\tau_n,L_{x_{(n+1)}}]\tau_n^{-1}
     \left(S^{-1}(y_{(-1)})h^0\otimes S(x_{(n)})h^1\otimes
           \cdots\otimes S(x_{(1)})h^n\otimes y_{(0)}\right)
     \label{fundamental}\\
 = & \left(S^{-1}(y_{(-1)})xh^0\otimes h^1\otimes\cdots\otimes
   h^n\otimes y_{(0)}\right)\nonumber\\
   & \hspace{.5cm} - \left(x_{(n+1)}S^{-1}(y_{(-1)})h^0\otimes
   x_{(n+2)}S(x_{(n)})h^1\otimes \cdots\otimes
   x_{(2n+1)}S(x_{(1)})h^n\otimes y_{(0)}\right)
\end{align*}
for any $(h^0\otimes\cdots\otimes h^n\otimes y)$ from $\B{T}_n(H,Y)$,
for any $x\in H$ and for any $n\geq 0$.
\end{defn}

\begin{defn}
An $H$--comodule $Y$ is called cocommutative iff 
\begin{align*}
(y_{(0)}\otimes y_{(\sigma 1)}\otimes\cdots\otimes y_{(\sigma n)})
 = (y_{(0)}\otimes y_{(1)}\otimes\cdots\otimes y_{(n)})
\end{align*}
for any $y\in Y$, for any $n\geq 1$ and for any $\sigma\in\Sigma_n$.
and element $x\in H$ is said to be in the cocenter of $H$ iff 
\begin{align*}
(x_{(\sigma 1)}\otimes\cdots\otimes x_{(\sigma n)})
 = (x_{(1)}\otimes\cdots\otimes x_{(n)})
\end{align*}
for any $n\geq 1$ and for any $\sigma\in\Sigma_n$.
\end{defn}

\begin{lem}\label{freakish}
Assume $Y$ is a cocommutative $H$--comodule and let $x\in H$ be an
element from the cocenter.  Consider $Y$ as a $H$--module via the counit
$\epsilon$.  Then
\begin{align*}
\tau_{n+1}^j\kappa_x\tau_n^{-j}(h^0\otimes\cdots\otimes h^n
    \otimes y)
 = & \left(h^0\otimes\cdots\otimes [S^{-1}(y_{(-1)}),x]h^j
           \otimes\cdots\otimes h^n\otimes y_{(0)}\right)
\end{align*}
is sent to $0$ by the morphism $\B{T}_*(H,Y)\xra{q_*}\B{PCM}_*(H,Y)$ and
also by the morphism $\B{T}_*(H,Y)\xra{p_*}\B{CM}_*(H,Y)$ or any
$(h^0\otimes\cdots\otimes h^n \otimes y)$ from $\B{CM}_n(H,Y)$ and for
any $0\leq j\leq n$.
\end{lem}
\begin{proof}
Regardless of $Y$ being cocommutative for any $(h^0\otimes\cdots\otimes
h^n\otimes y)$ from $\B{T}_n(H,Y)$, for any $x$ from $H$, for any $n\geq
0$ and $0\leq j\leq n$ consider the expression
\begin{align}
\tau_n^j\kappa_x\tau_n^{-j} & (h^0\otimes\cdots\otimes h^n
    \otimes y)\nonumber\\
 = & \tau_n^j\kappa_x(h^j\otimes\cdots\otimes h^n\otimes y_{(-j)}h^0
                \otimes\cdots\otimes y_{(-1)}h^{j-1}\otimes 
                y_{(0)})\nonumber\\
 = & \tau_n^j[\tau_n,L_{x_{(n+1)}}]
     \left(S^{-1}(y_{(-1)})h^j\otimes S(x_{(n)})h^{j+1}
           \otimes\cdots\right.\nonumber\\
   & \hspace{2cm}\left.\otimes S(x_{(j+1)})h^n\otimes
           S(x_{(j)})y_{(-j-1)}h^0\otimes\cdots\otimes 
           S(x_{(1)})y_{(-2)}h^{j-1}\otimes y_{(0)}\right)\nonumber\\
 = & \tau_n^j\left(S^{-1}(y_{(-1)})xh^j\otimes h^{j+1}\otimes\cdots
     \otimes h^n\otimes y_{(-j-1)}h^0\otimes\cdots\otimes
             y_{(-2)}h^{j-1}\otimes y_{(0)}\right)\nonumber\\
   & \hspace{.5cm} 
     - \tau_n^j\left(x_{(n+1)}S^{-1}(y_{(-1)})h^j\otimes
             x_{(n+2)}S(x_{(n)})h^{j+1}\otimes \cdots\right.\nonumber\\
   & \hspace{1.5cm}\otimes x_{(2n+1-j)}S(x_{(j+1)})h^n\otimes 
             x_{(2n+2-j)}S(x_{(j)})y_{(-j-1)}h^0
             \otimes\cdots\\
   & \left.\hspace{1.5cm}\otimes x_{(2n+1)}S(x_{(1)})y_{(-2)}h^{j-1}
             \otimes y_{(0)}\right)\nonumber
\end{align}
which finally reduces to 
\begin{align}
\tau_n^j\kappa_x\tau_n^{-j} & (h^0\otimes\cdots\otimes h^n
    \otimes y)\nonumber\\
 = & \left(S^{-1}(y_{(-1)})y_{(-2j-1)}h^0\otimes\cdots\otimes
           S^{-1}(y_{(-j)})y_{(-j-2)}h^{j-1}\right.\nonumber\\
   & \left.\hspace{1.5cm}\otimes
           S^{-1}(y_{(-j-1)})xh^j\otimes h^{j+1}\otimes\cdots
           \otimes h^n\otimes y_{(0)}\right)\label{Kappa}\\
   & - \left(S^{-1}(y_{(-1)})x_{(2n+2-j)}S(x_{(j)})y_{(-2j-1)}h^0
           \otimes\cdots\otimes
           S^{-1}(y_{(-j)})x_{(2n+1)}S(x_{(1)})y_{(-j-2)}h^{j-1}
           \right.\nonumber\\
   & \left.\hspace{1.5cm}\otimes
           x_{(n+1)}S^{-1}(y_{(-j-1)})h^j\otimes 
           x_{(n+2)}S(x_{(n)})h^{j+1}\otimes\cdots
           \otimes x_{(2n+1-j)}S(x_{(j+1)})h^n\otimes 
           y_{(0)}\right)\nonumber
\end{align}
Now, one can see that if $Y$ is cocommutative $H$--comodule and $x$ is
in the cocenter of $H$ then
\begin{align}
\tau_n^j\kappa_x\tau_n^{-j} & (h^0\otimes\cdots\otimes h^n
    \otimes y)
 = \left(h^0\otimes\cdots\otimes [S^{-1}(y_{(-1)}),x]h^j
           \otimes\cdots\otimes h^n\otimes y_{(0)}\right)
\end{align}
as I wanted to show.  Observe that $p_*[\tau_n^j,L_x]\equiv 0$ for any
$x\in H$ and therefore $p_*\kappa_x\equiv 0$.  This means
$p_n\tau_n^j\kappa_x\tau_n^{-j} = t_n^jp_n\kappa_x\tau_n^j=0$ for any
$0\leq j\leq n$ and $x$ from the cocenter of $H$.
\end{proof}

\section{Computations}\label{Computations}

\subsection{$\C{H}(n)$: Transverse geometry in codimension $n$}
\label{TransverseCodimensionN}

The Hopf algebra $\C{H}(n)$ as an algebra is generated by the elements
$X_k, Y_i^j$ and $\delta_{bc; i_1,\ldots,i_m}^a$ subject to the
following relations (I use Einstein summation notation)
\begin{align}
[X_k,X_\ell] = & R^i_{jk\ell} Y_i^j\label{LieAlg1}\\
[X_k,Y_i^j]  = & \delta^j_k X_i\\
[Y_i^j, Y_k^\ell]     = & \delta_k^j Y_i^\ell - \delta_j^\ell Y^i_k\label{LieAlg2}\\
[Y_i^j,\delta_{bc}^a] = & \delta_b^j\delta_{ic}^a -
                          \delta_i^a\delta_{bc}^j\label{LieAlg3}\\
[X_{i_m},\ldots,[X_{i_1},\delta_{bc}^a]\ldots]
                      = & \delta_{bc;i_1,\ldots,i_m}^a\label{Hn_inductive}
\end{align}
where $\delta_\alpha^\beta$ are the Kr\"onecker's $\delta$ functions.
The comultiplication on the generators $X_k$, $Y_i^j$ and
$\delta_{bc}^a$ are defined as
\begin{align}
\Delta(Y_i^j) = & (Y_i^j\otimes 1) + (1\otimes Y_i^j)\\
\Delta(\delta_{bc}^a) = & (\delta_{bc}^a\otimes 1) +
                          (1\otimes\delta_{bc}^a)\\
\Delta(X_i)   = & (X_i\otimes 1) + (1\otimes X_i) 
                  + (\delta_{ij}^k\otimes Y_k^j)
\end{align}
Comultiplication on $\delta_{bc;i_1,\ldots,i_m}^a$ is defined
inductively by using (\ref{Hn_inductive}).  This presentation is taken
from \cite{ConnesMoscovici:HopfCyclicCohomology}

Let $U(\G{gl}_n)$ be the universal enveloping algebra of the Lie algebra
generated by the symbols $\{Y_i^j\}_{i,j}$ subject to the condition
stated in Equation~(\ref{LieAlg2}).  Also let $\C{D}(n)$ be the
polynomial algebra on the symbols $\{\delta_{bc;I}^a\}_{a,b,c,I}$.  It
is clear that both $U(\G{gl}_n)$ and $\C{D}(n)$ are sub-coalgebras of
$\C{H}(n)$.

\subsubsection{Coefficients in $\C{D}(n)$}

Let $\C{D}_m(n)$ be the sub-Hopf algebra of $\C{D}(n)$ generated by
symbols of the form $\delta_{bc;I}^a$ where $|I|\leq m$.  It is clear
that $\C{D}_m(n)$ is a subalgebra, but not so clear that it is a
sub-coalgebra.  We prove this by induction on the length of the
multi-indices $I$ in $\delta_{bc;I}^a$.  For $m=0$, since
$\delta_{bc}^a$ is primitive the statement easily follows.  So, assume
for $|I|<m$ the statement holds.  Take $\delta_{bc;Ii}^a =
[X_i,\delta_{bc;I}^a]$ from $\C{D}_{m}(n)$.  Then
\begin{align*}
\Delta(\delta_{bc;Ii}^a)
 = & \Delta[X_i,\delta_{bc;I}^a]\\
 = & \left(X_{i,(1)}\delta_{bc;I,(1)}^a\otimes 
           X_{i,(2)}\delta_{bc;I,(2)}^a\right)
     - \left(\delta_{bc;I,(1)}^aX_{i,(1)}\otimes
             \delta_{bc;I,(2)}^aX_{i,(2)}\right)\\
 = & \left(\delta_{bc;I,(1)}^a\otimes X_i\delta_{bc;I,(2)}^a\right)
     + \left(X_i\delta_{bc;I,(1)}^a\otimes\delta_{bc;I,(2)}^a\right)
     + \left(\delta_{ij}^k\delta_{bc;I,(1)}^a\otimes
             Y_k^j\delta_{bc;I,(2)}^a\right)\\
   & - \left(\delta_{bc;I,(1)}^a\otimes \delta_{bc;I,(2)}^aX_i\right)
     - \left(\delta_{bc;I,(1)}^aX_i\otimes\delta_{bc;I,(2)}^a\right)
     - \left(\delta_{bc;I,(1)}^a \delta_{ij}^k\otimes
             \delta_{bc;I,(2)}^aY_k^j\right)\\
 = & \left(\delta_{bc;I,(1)}^a\otimes [X_i,\delta_{bc;I,(2)}^a]\right)
     + \left([X_i,\delta_{bc;I,(1)}^a]\otimes\delta_{bc;I,(2)}^a\right)
     + \left(\delta_{ij}^k\delta_{bc;I,(1)}^a\otimes
             [Y_k^j,\delta_{bc;I,(2)}^a]\right)
\end{align*}
Taking commutators of $\delta_{bc;I,(\alpha)}$ with $X_i$ raise the
length of the multi index by one.  However, the induction hypothesis
tells us the length of the multi index in $\delta_{bc;I,(\alpha)}$ is at
most $m-1$.  Therefore the length of the multi-index in
$[X_i,\delta_{bc;I,(\alpha)}^a]$ is at most $m$.  Since taking
commutators of $\delta_{bc;I,(\alpha)}$ with $Y_j^i$ does not change the
length of the multi-indices, the result follows.

Since $\C{D}_0(n)$ is cocommutative coalgebra of $\C{H}(n)$ and $Y_i^j$
is in the cocenter of $\C{H}(n)$, I have
\begin{align*}
\tau_{n+1}^s\kappa_{Y_i^b}\tau_n^{-s}
     (h^0\otimes\cdots\otimes h^n\otimes \delta_{bc}^a)
 = & (h^0\otimes\cdots\otimes [Y_i^b,\delta_{bc}^a]h^s\otimes
      \cdots\otimes h^n\otimes 1)\\
 = & (h^0\otimes\cdots\otimes \delta_{ic}^a h^s\otimes\cdots\otimes
      h^n\otimes 1)
\end{align*}
Now, assume terms of the form 
\begin{align*}
(h^0\otimes\cdots\otimes \delta_{ic}^a h^s\otimes\cdots\otimes
      h^n\otimes \delta_I)
\end{align*}
are sent to $0$ under $p_*$ for any $|I|<m$.  Take
$(h^0\otimes\cdots\otimes h^n\otimes \delta_{bc}^a\delta_A)$ with
$|A|=m$ and consider
\begin{align*}
\tau_{n+1}^s\kappa_{Y_i^b}\tau_n^{-s}
   & (h^0\otimes\cdots\otimes h^n\otimes \delta_{bc}^a\delta_A)\\
 = & (h^0\otimes\cdots\otimes [Y_i^b,\delta_{bc}^a]h^s\otimes
      \cdots\otimes h^n\otimes \delta_A)
    + \sum_{I\subseteq A}
     (h^0\otimes\cdots\otimes [Y_i^b,\delta_{bc}^a\delta_I]h^s\otimes
      \cdots\otimes h^n\otimes \delta_{A\setminus I})
\end{align*}
are all sent to zero.  Note that since
\begin{align*}
[Y_i^j,\delta_{b_1c_1}^{a_1}\cdots\delta_{b_wc_w}^{a_w}]
 = & [Y_i^j,\delta_{b_1c_1}^{a_1}]
       \delta_{b_2c_2}^{a_2}\cdots\delta_{b_wc_w}^{a_w}
     + \delta_{b_1c_1}^{a_1}
       [Y_i^j,\delta_{b_2c_2}^{a_2}\cdots\delta_{b_wc_w}^{a_w}]
\end{align*}
second summand is sent to $0$ under $p_*$ by induction hypothesis.  Then
so is the first term, as I wanted to show.  Therefore, by the help of
Lemma~\ref{quotient}, I can use the $\C{H}(n)$--module coalgebra
\begin{align*}
\C{H}(n)/\left<\delta_{bc}^a\right|_{a,b,c}
\end{align*}

On this quotient $X_i$'s act like primitive elements.  Thus elements of
the form 
\begin{align*}
\tau_{n+1}^s\kappa_{X_i}\tau_n^{-s}
     (h^0\otimes\cdots\otimes h^n\otimes \delta_{bc}^a)
 = & (h^0\otimes\cdots\otimes [X_i,\delta_{bc}^a]h^s\otimes
      \cdots\otimes h^n\otimes 1)\\
 = & (h^0\otimes\cdots\otimes \delta_{bc;i}^a h^s\otimes\cdots\otimes
      h^n\otimes 1)
\end{align*}
are sent to $0$ by $p_*$.  Then, assume by induction that terms of the
form
\begin{align*}
(h^0\otimes\cdots\otimes \delta_{bc;i}^ah^s\otimes
        \cdots\otimes h^n\otimes \delta_I)
\end{align*}
are sent to $0$ with $|I|<m$.  Then for $|A|=m$
\begin{align*}
\tau_{n+1}^s\kappa_{X_i}\tau_n^{-s}
   &  (h^0\otimes\cdots\otimes h^n\otimes \delta_{bc}^a\delta_A)\\
 = & (h^0\otimes\cdots\otimes [X_i,\delta_{bc}^a]h^s\otimes
      \cdots\otimes h^n\otimes \delta_A)
     + \sum_{I\subseteq A}
       (h^0\otimes\cdots\otimes [X_i,\delta_{bc}^a\delta_I]h^s\otimes
        \cdots\otimes h^n\otimes \delta_{A\setminus I})
\end{align*}
is sent to zero.  Second summand is sent to $0$ by induction hypothesis
and the argument I presented above.  This means, instead of $\C{H}(n)$,
I can use
\begin{align*}
\C{H}(n)/\left<\delta_{bc}^a,\delta_{jk;\ell}^i\right|_{a,b,c,i,j,k,\ell}
\end{align*}
Now, by the little Lemma I proved at the beginning of this section, the
elements in $\C{D}_1(n)$ behave like primitive elements.  Then by
induction, instead of $\C{H}(n)$, I can use
\begin{align}
\C{H}(n)/\left<\delta_{jk;I}^i\right|_{a,b,c,i,j,k,I}
\end{align}
Since $\C{H}(n)$ has a basis of the form $\{\delta_IX_JY_K\}_{I,J,K}$
the quotient is isomorphic to the universal enveloping algebra of the
Lie algebra $\G{a}_n$ generated by $\{X_i, Y_j^k\}_{i,j,k}$ subject to
relations given in Equation~(\ref{LieAlg1}) to Equation~(\ref{LieAlg3}).
Then
\begin{align}
\B{CM}_*(\C{H}(n),\C{D}(n))
\cong {}_{\C{H}(n)}\B{T}_*(U(\G{a}_n),\C{D}(n))
\cong {}_{U(\G{a}_n)}\B{T}_*(U(\G{a}_n),k)\otimes\C{D}(n)
\cong \B{CM}_*(U(\G{a}_n),k)\otimes\C{D}(n)
\end{align}
Therefore, by using \cite{ConnesMoscovici:HopfCyclicCohomology} or
\cite{Crainic:CyclicCohomologyOfHopfAlgebras} we obtain
\begin{align}
HP^\B{CM}_n(\C{H}(n),\C{D}(n)) := HP_n \B{CM}_*(\C{H}(n),\C{D}(n))
\cong \bigoplus_{i \equiv n \text{ mod } 2}H^{Lie}_i(\G{a}_n,k)\otimes\C{D}(n)
\end{align}

\subsubsection{Coefficients in $U(\G{gl}_n)$}

$U(\G{gl}_n)$ is a cocommutative sub-coalgebra of $\C{H}(n)$.  We also
consider $U(\G{gl}_n)$ as a $\C{H}(n)$--module through $\epsilon$.  Then
by Lemma~\ref{freakish} any element of the form
\begin{align*}
\tau_{n+1}^s\kappa_{Y_a^i}\tau_n^{-s}
     (h^0\otimes\cdots\otimes h^n\otimes Y_i^j)
 = & (h^0\otimes\cdots\otimes [Y_a^i, Y_i^j]h^s\otimes\cdots\otimes
      h^n\otimes 1)\\
 = & (h^0\otimes\cdots\otimes Y_a^j h^s\otimes\cdots\otimes
      h^n\otimes 1)
\end{align*}
is sent to zero under $p_*$ for any $0\leq s\leq n$ and
$a,i,j,k=1,\ldots, n$.  Assume (as an induction hypothesis) that terms
of the form
\begin{align*}
(h^0\otimes\cdots\otimes Y_i^j h^s\otimes\cdots\otimes
      h^n\otimes Y_I)
\end{align*}
are sent to $0$ by $p_*$ where length of the multi-index $|I|<m$.  Take
$(h^0\otimes\cdots\otimes Y_a^j h^s\otimes\cdots\otimes h^n\otimes
Y_i^jY_A)$ with $|A|=m$.  Then
\begin{align*}
\tau_{n+1}^s\kappa_{Y_a^i}\tau_n^{-s}
   & (h^0\otimes\cdots\otimes h^n\otimes Y_i^jY_A)\\
 = & (h^0\otimes\cdots\otimes [Y_a^i, Y_i^j]h^s\otimes\cdots\otimes
      h^n\otimes Y_A)
     + \sum_{I\subseteq A}
       \left(h^0\otimes\cdots\otimes [Y_a^i, Y_i^jY_I]h^s\otimes\cdots
           \otimes h^n\otimes Y_{A\setminus I}\right)
\end{align*}
are all sent to $0$ where the sum is taken over all {\em ordered}
subsets of $A$.  The second sum is sent to $0$ by induction hypothesis,
then so is the first term.  The result follows. Then instead of
$\C{H}(n)$, I can use
\begin{align}
\C{H}(n)/\left<Y_i^j\right|_{i,j}
\end{align}
where $\left<z_\lambda\right|_{\lambda\in\Lambda}$ denotes the left
ideal in $H$ generated by the set $\{z_\lambda|\ \lambda\in\Lambda\}$.

Now consider
\begin{align*}
\tau_{n+1}^s\kappa_{\delta_{bc}^a}\tau_n^{-s}
     (h^0\otimes\cdots\otimes h^n\otimes Y_a^j)
 = & (h^0\otimes\cdots\otimes [\delta_{bc}^a, Y_a^j]h^s\otimes
      \cdots\otimes h^n\otimes 1)\\
 = & (h^0\otimes\cdots\otimes \delta_{bc}^j h^s\otimes\cdots\otimes
      h^n\otimes 1)
\end{align*}
which is sent to $0$ by $p_*$.  By a similar argument I presented above,
one can conclude that
\begin{align*}
(h^0\otimes\cdots\otimes \delta_{bc}^j h^s\otimes\cdots\otimes
      h^n\otimes Y_A)
\end{align*}
are sent to $0$, for any $0\leq s\leq n$, $b,c,j=1,\ldots,n$ and
multi-index $A$.  Thus one can further reduce $\C{H}(n)$ to
\begin{align}
\C{H}(n)/\left<Y_i^j, \delta_{bc}^a\right|_{i,j,a,b,c}
\end{align}
Now, on this quotient, elements of the form $\delta_{bc;d}^a$ act like
primitive elements.  Consider
\begin{align*}
\tau_{n+1}^s\kappa_{\delta_{bc;d}^a}\tau_n^{-s}
    & (h^0\otimes\cdots\otimes h^n\otimes Y_i^j)\\
  = & (h^0\otimes\cdots\otimes [\delta_{bc;d}^a, Y_i^j]h^s\otimes
       \cdots\otimes h^n\otimes 1)\\
  = & (h^0\otimes\cdots\otimes [[X_d, \delta_{bc}^a],Y_i^j]h^s\otimes
      \cdots\otimes h^n\otimes 1)\\
  = & (h^0\otimes\cdots\otimes [X_d,[\delta_{bc}^a,Y_i^j]]h^s\otimes
      \cdots\otimes h^n\otimes 1)
    - (h^0\otimes\cdots\otimes [\delta_{bc}^a,[X_d,Y_i^j]]h^s\otimes
      \cdots\otimes h^n\otimes 1)\\
  = & -\delta_b^j(h^0\otimes\cdots\otimes\delta_{ic;d}^a h^s\otimes
                  \cdots\otimes h^n\otimes 1)
      +\delta_a^i(h^0\otimes\cdots\otimes\delta_{bc;d}^j h^s\otimes
                  \cdots\otimes h^n\otimes 1)\\
    & + \delta_d^j(h^0\otimes\cdots\otimes\delta_{bc;i}^a h^s\otimes
         \cdots\otimes h^n\otimes 1)
\end{align*}
Since all indices are arbitrary, one can see that, elements of the form
\begin{align*}
(h^0\otimes\cdots\otimes\delta_{bc;d}^a h^s\otimes\cdots\otimes h^n\otimes 1)
\end{align*}
are sent to $0$ by $p_*$.  By induction on both $A$ and $I$ one can show
that elements of the form 
\begin{align*}
(h^0\otimes\cdots\otimes\delta_{bc;I}^a h^s\otimes\cdots\otimes h^n\otimes Y_A)
\end{align*}
are all sent to $0$ by $p_*$.  Then one can use the $\C{H}(n)$--module
coalgebra
\begin{align*}
\C{H}(n)/\left<Y_i^j, \delta_{bc;I}^a\right|_{i,j,a,b,c,I}
\end{align*}

On this quotient, $X_i$'s act like primitive elements.  Thus
\begin{align*}
\tau_{n+1}^s\kappa_{X_j}\tau_n^{-s}
     (h^0\otimes\cdots\otimes h^n\otimes Y_i^j)
 = & (h^0\otimes\cdots\otimes [X_j, Y_i^j]h^s\otimes
      \cdots\otimes h^n\otimes 1)\\
 = & (h^0\otimes\cdots\otimes X_i h^s\otimes\cdots\otimes
      h^n\otimes 1)
\end{align*}
which is sent to $0$ by $p_*$.  Again, induction on the length of
multi-indices allows us to conclude that
\begin{align*}
(h^0\otimes\cdots\otimes X_i h^s\otimes\cdots\otimes
      h^n\otimes Y_A)
\end{align*}
is sent to $0$ for any $0\leq s\leq n$, $i=1,\ldots,n$ and multi-index
$A$.  Therefore, finally, $\C{H}(n)$ reduces to
\begin{align}
\C{H}(n)/\left<X_k, Y_i^j, \delta_{bc;I}^a\right|_{i,j,k,a,b,c,I}
\cong k
\end{align}
which implies
\begin{align}
\B{CM}_*(\C{H}(n),U(\G{gl}_n))
\cong {}_{\C{H}(n)}\B{T}_*(k,U(\G{gl}_n))
\cong \B{T}_*(k,k)\otimes U(\G{gl}_n)
\cong CC_*(k)\otimes U(\G{gl}_n)
\end{align}
Therefore
\begin{align}
HC^\B{CM}_*(\C{H}(n),U(\G{gl}_n)) := HC_*\B{CM}_*(\C{H}(n),U(\G{gl}_n))
\cong HC_*(k)\otimes U(\G{gl}_n)
\end{align}

\subsection{Quantum deformations of universal enveloping algebras}\label{quantum}

Let $\G{g}$ be a finite dimensional (semisimple) Lie algebra over
$k=\B{C}$.  Fix a Cartan subalgebra $\C{H}$ and let $(a_{ij})$ be the
corresponding Cartan matrix.  Then, the quantum deformation $U_q\G{g}$
is generated by the elements $\{K_i^\pm, X_i^\pm|\ i=1,\ldots,N-1\}$
subject to certain relations.  Among those, the following relations are
of importance:
\begin{align}
K_i^+K_i^-          = & K_i^-K_i^+ = 1\label{Ug1}\\
K_i^+ X_j^\pm K_i^- = & q^{\pm a_{ij}} X_j^\pm\label{Ug2}\\
[X_i^+,X_j^-]       = & \frac{\delta_i^j}{q-q^{-1}}(K_i^+-K_i^-)\label{Ug3}\\
[K_i,K_j]           = & 0
\end{align}
The comultiplication is defined on generators as
\begin{align}
\Delta(K_i)   = & (K_i^+\otimes K_i^+)\\
\Delta(X_i^+) = & (1\otimes X_i^+) + (X_i^+\otimes K_i^+)\\
\Delta(X_i^-) = & (K_i^-\otimes X_i^-) + (X_i^-\otimes 1)
\end{align}
for any $i=1,\ldots,N-1$.  From these definitions one can deduce that
\begin{align*}
\epsilon(X^\pm_i) = & 0 &
\epsilon(K^\pm_i) = & 1 &
S(K^\pm_i)        = & K^\mp_i &
S(X^+_i)          = & -X^+_iK^-_i&
S(X^-_i)          = & -K^+_iX^-_i
\end{align*}
Considering Equation~(\ref{Kappa}) in Lemma~\ref{freakish} with $j=0$, I
see that terms of the form
\begin{align}
\kappa_x(h^0\otimes\cdots\otimes h^n\otimes y)
 = & \left(S^{-1}(y_{(-1)})xh^0\otimes h^1\otimes\cdots\otimes
   h^n\otimes y_{(0)}\right)\\ 
   & -  \left(x_{(n+1)}S^{-1}(y_{(-1)})h^0\otimes
              x_{(n+2)}S(x_{(n)})h^1\otimes \cdots\otimes
              x_{(2n+1)}S(x_{(1)})h^n\otimes y_{(0)}\right)\label{killerE}
\end{align}
is mapped to zero under $p_*$.  For $x=K^\pm_i$, this difference is
\begin{align*}
\left(S^{-1}(y_{(-1)})K^\pm_i h^0\otimes
           h^1\otimes\cdots\otimes h^n\otimes  y_{(0)}\right)
 - \left(K^\pm_i S^{-1}(y_{(-1)})h^0\otimes
           h^1\otimes\cdots\otimes h^n\otimes  y_{(0)}\right)
\end{align*}
For $X^+_i$, I have
\begin{align*}
\Delta^{2n}(X^+_i)
 = & \sum_{j=0}^{2n}(1^{\otimes 2n-j}\otimes X^+_i\otimes 
                     (K^+_i)^{\otimes j})
\end{align*}
Then, for $x=X^+_i$
\begin{align*}
& \left(x_{(n+1)}S^{-1}(y_{(-1)})h^0\otimes x_{(n+2)}S(x_{(n)})\otimes
           \cdots\otimes x_{(2n+1)}S(x_{(1)})h^n\otimes y_{(0)}\right)\\
& \hspace{1cm}  =  L_{X^+_i}\left(S^{-1}(y_{(-1)})
                    h^0\otimes\cdots\otimes h^n\otimes y_{(0)}\right)\\
& \hspace{1.5cm} - \sum_{j=1}^n q^2\left(K^+_iS^{-1}(y_{(-1)})h^0
         \otimes h^1\otimes\cdots\otimes 
         h^{j-1}\otimes X^+_i h^j\otimes K^+_i h^{j+1}
         \otimes\cdots\otimes K^+_i h^n\otimes y_{(0)}\right)\\
& \hspace{1cm} =  L_{X^+_i}\left(S^{-1}(y_{(-1)})h^0\otimes
         \cdots\otimes h^n\otimes y_{(0)}\right)\\
& \hspace{1.5cm}  - q^{2 a_{ii}} L_{X^+_i}\left(K^+_i S^{-1}(y_{(-1)})h^0
            \otimes h^1\otimes\cdots\otimes h^n\otimes y_{(0)}\right)\\
& \hspace{1.5cm}  + q^{2 a_{ii}}(X^+_iK^+_i S^{-1}(y_{(-1)})h^0\otimes 
         h^1\otimes\cdots\otimes h^n\otimes y_{(0)})\\
& \hspace{1cm} =  L_{X^+_i}\left((1-q^{2 a_{ii}}K^+_i)S^{-1}(y_{(-1)})h^0
          \otimes\cdots\otimes h^n\otimes y_{(0)}\right)\\
& \hspace{1.5cm}  + \left(K^+_iX^+_i S^{-1}(y_{(-1)})
         h^0\otimes h^1\otimes\cdots
             \otimes h^n\otimes y_{(0)}\right)
\end{align*}
Therefore,
\begin{align*}
\kappa_{X^+_i} & (h^0\otimes\cdots\otimes h^n\otimes y)\\
 = & \left(S^{-1}(y_{(-1)})X^+_ih^0\otimes h^1\otimes\cdots\otimes
           h^n\otimes y_{(0)}\right)
    - \left(K^+_iX^+_i S^{-1}(y_{(-1)})h^0\otimes h^1\otimes\cdots\otimes 
            h^n\otimes y_{(0)}\right)\\
   & - L_{X^+_i}\left((1-q^{2 a_{ii}}K^+_i)S^{-1}(y_{(-1)})h^0\otimes\cdots
          \otimes h^n\otimes y_{(0)}\right)
\end{align*}
must be sent to $0$ in $\B{CM}_*(U_q(\G{g}), Y)$.  Also, instead of
using $X^+_i$, if I used $X^-_iK^+_i$, one can see that
\begin{align*}
\kappa_{X^-_iK^+_i} & (h^0\otimes\cdots\otimes h^n\otimes y)\\
 = & \left(S^{-1}(y_{(-1)})X^-_iK^+_ih^0\otimes h^1\otimes\cdots\otimes
           h^n\otimes y_{(0)}\right)
   - \left(K^+_iX^-_iK^+_i S^{-1}(y_{(-1)})h^0\otimes h^1\otimes\cdots\otimes 
            h^n\otimes y_{(0)}\right)\\
   & - L_{X^-_iK^+_i}\left((1-q^{-2a_{ii}}K^+_i)S^{-1}(y_{(-1)})h^0\otimes\cdots
          \otimes h^n\otimes y_{(0)}\right)
\end{align*}
must be sent to 0 in $\B{CM}_*(U_q(\G{g}),Y)$ under $p_*$. Since
$\epsilon(X^\pm_i)=0$, I must have
\begin{align}
p_*\left(S^{-1}(y_{(-1)})K^\pm_i h^0\otimes
           h^1\otimes\cdots\otimes h^n\otimes  y_{(0)}\right)
   = &  p_*\left(K^\pm_iS^{-1}(y_{(-1)}) h^0\otimes
           h^1\otimes\cdots\otimes h^n\otimes  y_{(0)}\right)
          \label{fundamental1}\\
p_*\left(S^{-1}(y_{(-1)})X^\pm_i h^0\otimes h^1\otimes\cdots\otimes h^n
      \otimes y_{(0)}\right)
   = & p_*\left(K^+_iX^\pm_i S^{-1}(y_{(-1)})h^0\otimes h^1\otimes
          \cdots\otimes h^n\otimes y_{(0)}\right)\label{fundamental2}
\end{align}
in $\B{CM}_*(U_q(\G{g}), Y)$. 

\subsubsection{Coefficients in $U_qsl(2)$}

Fix $1\leq j\leq N-1$ and let $Y$ be to the copy of $U_qsl(2)$ generated
by the symbols $X^\pm_j$ and $K^\pm_j$ subject to the conditions stated
in Equations~(\ref{Ug1}) to (\ref{Ug3}).

Use $y=X^+_j$ and the fact that $S^{-1}(X^+_j)= - K^-_j X^+_j$ in
Equation~(\ref{fundamental1}) to obtain
\begin{align*}
p_*\left(K^-_j X^+_jK^\pm_i h^0\otimes
           h^1\otimes\cdots\otimes h^n\otimes  K^+_j\right)
   = &  p_*\left(K^\pm_i K^-_j X^+_j h^0\otimes
           h^1\otimes\cdots\otimes h^n\otimes  K^+_j\right)
\end{align*}
By a clever choice of $(h^0\otimes\cdots\otimes h^n)=(K^\mp_i u^0\otimes
K^+_j h^1\otimes\cdots\otimes K^+_j u^n)$, I can further reduce the
equality above to
\begin{align*}
K^+_j\cdot p_*&\left(K^+_jX^+_ju^0\otimes
           K^+_ju^1\otimes\cdots\otimes K^+_ju^n\otimes  K^+_j\right)\\
  = & K^+_j\cdot q^{\pm a_{ij}}p_*\left(K^+_jX^+_ju^0\otimes
           K^+_ju^1\otimes\cdots\otimes K^+_ju^n\otimes  K^+_j\right)
\end{align*}
for any $i$ and $s$.  Since $i$ is arbitrary and there is at least $i$
for which $a_{ij}\neq 0$, I conclude that
\begin{align*}
p_*\left(X^+_j u^0\otimes u^1\otimes\cdots\otimes u^n\otimes 1\right)
 = & 0
\end{align*}
By using $y=K^+_jX^-_j$, one can get a similar result as above.
Therefore, I can safely say
\begin{align*}
p_*\left(u^0\otimes\cdots\otimes X^\pm_ju^s\otimes
           \cdots\otimes u^n\otimes  K^c_j\right)
 = & 0
\end{align*}
for any $0\leq s\leq n$ and $c\in\B{Z}$.  

Assume, by induction, that
\begin{align*}
p_*(X^+_jh^0\otimes h^1\otimes\cdots\otimes h^n\otimes K^c_j X^{+m}_j)
 = 0
\end{align*}
for any $c$ and $m<n$.  Consider Equation~(\ref{fundamental1}) with
$y=X^{+(n+1)}_j$ to conclude that
\begin{align*}
\sum_{s=1}^{n+1} &c_s(K^{-s}_jX^{+s}_jK^\pm_i h^0 \otimes h^1\otimes
                \cdots\otimes h^n\otimes K^{+s}_j X^{+(n+1-s)}_j)\\
   & \hspace{1.5cm}- c_s(K^\pm_iK^{-s}_jX^{+s}_j h^0\otimes h^1\otimes
                \cdots\otimes h^n\otimes K^{+s}_j X^{+(n+1-s)}_j)\\
 = & \sum_{s=1}^{n+1} (1-q^{s\cdot a_{ij}})c_s
          (K^{-s}_jX^{+s}_jK^\pm_i h^0 \otimes h^1\otimes
                \cdots\otimes h^n\otimes K^{+s}_j X^{+(n+1-s)}_j)
\end{align*}
is sent to $0$ by $p_*$ where $c_s>0$ is a constant determined by the
relations $K^+_j X^+j K^-_j=q^{a_{jj}} X^+_j$.  The induction hypothesis
leaves only one term
\begin{align*}
(1-q^{s\cdot a_{ij}})c_1
          (K^-_jX^+_jK^\pm_i h^0 \otimes h^1\otimes
                \cdots\otimes h^n\otimes K^+_j X^{+n}_j)
\end{align*}
which is sent to $0$.  Since $q$ is not a root of unity, $(1-q^{s\cdot
a_{ij}})\neq 0$.  Therefore one can conclude
\begin{align*}
p_*(X^+_jh^0 \otimes h^1\otimes
                \cdots\otimes h^n\otimes K^+_j X^{+n}_j)
 = & 0
\end{align*}
for any $n\geq 0$ as I wanted to show.

One can similarly show that 
\begin{align}\label{base}
p_*(X^\pm_jh^0 \otimes h^1\otimes
                \cdots\otimes h^n\otimes K^{\ell}_j X^{+n}_j X^{-m}_j)
 = & 0
\end{align}
for all $m,n\in\B{N}$ and $\ell\in\B{Z}$.  Now apply $t_*$ the the
equation above to obtain
\begin{align*}
p_*(S^{-1}(y_{(-1)})h^n\otimes X^\pm_j h^0\otimes\cdots\otimes h^{n-1}
    \otimes y_{(0)})
 = & 0
\end{align*}
where $y=K^\ell_jX^{+n}_jX^{-m}_j$.  In $S^{-1}(y_{(-1)})$ I can not
have $X^\pm_j$ appearing because of Equation~(\ref{base}).  Therefore
\begin{align*}
p_*(K^\ell_j h^n\otimes X^\pm_jh^0 \otimes h^1\otimes
      \cdots\otimes h^{n-1}\otimes K^{\ell}_j X^{+n}_j X^{-m}_j)
 = & 0
\end{align*}
which implies
\begin{align*}
p_*(u^0\otimes X^\pm_ju^1 \otimes u^2\otimes
      \cdots\otimes u^{n-1}\otimes K^{\ell}_j X^{+n}_j X^{-m}_j)
 = & 0
\end{align*}
for all $(u^0\otimes\cdots\otimes u^n)$ from $U_q(\G{g})^{\otimes n}$.
Again by induction
\begin{align*}
p_*(u^0\otimes\cdots\otimes X^\pm_ju^s\otimes
     \cdots\otimes u^n\otimes K^{\ell}_j X^{+n}_j X^{-m}_j)
 = & 0
\end{align*}
for any $s$ and for any $(u^0\otimes\cdots\otimes u^n)$ from
$U_q(\G{g})^{\otimes n}$.  Therefore, I can use the quotient
\begin{align}
U:= 
\ U_q(\G{g})\slash\left<X^\pm_j\right|
\ \cong \ k[C_2]\ltimes U_q(\G{g}_j)
\end{align}
where $U_q(\G{g}_j)$ is obtained from $U_q(\G{g})$ by deleting all
occurrences of $\{K^{\pm m}_j, X^{\pm n}_j\}_{m,n}$.  The crossed product
structure $k[C_2]\ltimes U_q(\G{g}_j)$ comes from the relation
$K^+_jX^\pm_iK^-_j=q^{\pm a_{ji}}X_i$ where $C_2=\left<K^\pm_j|\ K^{\pm
2}_j\right>$.

Let $Y_0=k\{K^\ell_j X^{+n}_jX^{-m}_j|\ \ell\text{ is even}\}$ and
$Y_1=k\{K^\ell_j X^{+n}_jX^{-m}_j|\ \ell\text{ is odd}\}=
K^+_jY_0$. Then by Lemma~\ref{quotient},
\begin{align}
\B{CM}_*(U_q(\G{g}),Y) 
\cong &  {}_{U_q(\G{g})}\B{T}_*(U,Y_0)\oplus 
         {}_{U_q(\G{g})}\B{T}_*(U,K^+_jY_0)\\
\cong &  {}_U\B{T}_*(U,k)\otimes Y_0\oplus 
         {}_U\B{T}_*(U,k_{K^+_j})\otimes Y_0\\
\cong & \B{CM}_*\left(k[C_2]\ltimes U_q(\G{g}_j),k[C_2]\right)\otimes Y_0
\end{align}
which implies
\begin{align}
HC^\B{CM}_*(U_q(\G{g}),Y) := HC_*(U_q(\G{g}),Y)
\cong  HC^\B{CM}_*\left(k[C_2]\ltimes U_q(\G{g}_j),k[C_2]\right)\otimes Y_0
\end{align}

\vspace{2cm}
\noindent{\small\sc 
Department of Mathematics, Ohio State University,
Columbus, Ohio 43210, USA}

\noindent{\it E-mail address:}\ {\tt kaygun@math.ohio-state.edu}

\end{document}